\documentclass[10pt,reqno]{amsart}
\usepackage{amsfonts}
\usepackage{amsthm}
\usepackage{amsmath}
\usepackage{amssymb}
\usepackage[utf8]{inputenc}
\usepackage{enumerate}
\usepackage{hyperref}
\usepackage{thmtools}

\newtheoremstyle{mystyle}
  {}
  {}
  {\itshape}
  {}
  {\bfseries}
  {.}
  { }
  {\thmname{#1}\thmnumber{ #2}\thmnote{ (#3)}}

\theoremstyle{mystyle}
\newtheorem{Thm}{Theorem}[section]
\newtheorem{Lem}[Thm]{Lemma}
\newtheorem{Cor}[Thm]{Corollary}
\newtheorem{Prop}[Thm]{Proposition}
\newtheorem{Conj}[Thm]{Conjecture}

\theoremstyle{definition}
\newtheorem{Def}[Thm]{Definition}
\newtheorem{Ex}[Thm]{Example}

\theoremstyle{remark}
\newtheorem{Rmk}[Thm]{Remark}

\declaretheoremstyle[%
  spaceabove=3pt,
  spacebelow=10pt,
  headfont=\normalfont\itshape,%
  postheadspace=.5em,%
  qed=\qedsymbol%
]{mystyle2} 
\declaretheorem[name={Proof},style=mystyle2,unnumbered,
]{pf}

\let\emptyset\varnothing

\let\epsilon\varepsilon
\let\star\nothing
\let\c\nothing

\newcommand{\R}{\mathbb{R}}

\newcommand{\C}{\mathbb{C}}
\renewcommand{\S}{\mathbb{S}}
\newcommand{\T}{\mathcal{T}}
\newcommand{\I}{\mathcal{I}}
\newcommand{\diam}{\mathrm{diam}}

\newcommand{\vol}{\mathrm{vol}}

\newcommand{\M}{\mathcal{M}}
\newcommand{\c}{\mathrm{Cap}}
\newcommand{\star}{\mathrm{Star}}
\newcommand{\cone}{\mathrm{Cone}}
\newcommand{\A}{\mathcal{A}}
\newcommand{\id}{\mathrm{id}}
\newcommand{\Z}{\mathbb{Z}}

\title{On the union of essentially distinct $\delta$-tubes}

\author{Qiuyu Ren}
\address{Department of Mathematics, Massachusetts Institute of Technology, Cambridge, MA 02139, USA}
\email{rqy18@mit.edu}
\begin{document}

\maketitle

\begin{abstract}
    We say two $\delta$-tubes (dimension $\delta\times\cdots\times\delta\times1$) in $\R^n$ are essentially distinct if the measure of their intersection is smaller than a half of a single $\delta$-tube. For a collection of essentially distinct $\delta$-tubes, we give the asymptotically sharp lower bound for the measure of their union. Then we characterize all sharp examples. We will give a new measurement of convexity based on the X-ray transform.
\end{abstract}

\section{Introduction}
The Kakeya conjecture has long been a famous unsolved problem in analysis. We call a compact set $E$ in $\R^n$ a Kakeya set if it contains a translation of every unit segment in $\R^n$. For $n\ge2$, there exists Kakeya set of measure zero. However, the Kakeya conjecture states that:

\begin{Conj}[Kakeya set conjecture]
Every Kakeya set in $\R^n$ has Hausdorff and Minkowski dimension $n$.
\end{Conj}

Of course, the result for Hausdorff dimension implies the one for Minkowski dimension.

The conjecture is trivial for $n=1$, easy for $n=2$ but remains unsolved for any $n\ge3$. One can ask various different questions which are closely related to Kakeya conjecture. Most of these ``Kakeya family'' problems remain open.

Our work is closely to one of these conjectures, namely the Kakeya tube conjecture. To formulate it, let us first fix some notations.

Let $a\in\R^n$, $e\in\S^{n-1}/\{\pm1\}$. A \textbf{$\delta$-tube} in $\R^n$ with \textbf{center} $a$ and \textbf{direction} $e$ is defined to be the set 
\begin{equation}\label{tube}
T=T_{e}^{\delta}(a)=\{\,x\in\R^n\colon |(x-a)\cdot\hat e|\le1/2,\ |(x-a)-((x-a)\cdot\hat e)\hat e|\le\delta/2\,\},    
\end{equation}
where $\hat e$ is any preimage of $e$ in $\S^{n-1}$.

The \textbf{angle} between two directions $e_1,e_2$ (or two $\delta$-tubes with direction $e_1,e_2$) is defined to be $\arccos|\hat e_1\cdot\hat e_2|$.

We say a collection of tubes $\{T_i\}_{i\in I}$ (or a collection of directions $\{e_i\}_{i\in I}$) is \textbf{$\delta$-separated} if the angle between every two of them is greater than $\delta$. We say it is \textbf{maximal $\delta$-separated} if it is $\delta$-separated and $\#I\ge c\delta^{1-n}$ for some prescribed small constant $c>0$.

\begin{Conj}[Kakeya tube conjecture]
For any $\epsilon>0$, there exists $c_\epsilon>0$ such that for $\delta>0$ and any maximal $\delta$-separated collection $\{T_i\}_{i\in I}$ of $\delta$-tubes we have
\begin{equation}\label{e1}
\left|\bigcup_{i\in I}T_i\right|\ge c_\epsilon\delta^\epsilon.
\end{equation}
Where $|\cdot|$ denotes the Lebesgue measure.
\end{Conj}

The Minkowski version of the Kakeya set conjecture is a direct corollary of the Kakeya tube conjecture. More generally, if we replace the exponent $\epsilon$ in \eqref{e1} by $n-d+\epsilon$ for some $0\le d\le n$, the result would imply that every Kakeya set has Minkowski dimension at least $d$. Many partial progresses have been made. For example, Wolff \cite{wolff1995improved} showed this for $d=(n+2)/2$. Later, Katz and Tao \cite{katz2002new} showed this for $d=(4n+3)/7$. Still better results have been made, but we do not attempt to give a full list here.

We will not directly tackle this problem. Instead, we consider a variant of it. We say a collection of $\delta$-tubes $\{T_i\}_{i\in I}$ is \textbf{essentially distinct} if $|T_i\cap T_j|\le c_0|T_i|$ for every $i\ne j\in I$ where $0<c_0<1$ is a fixed constant. Notice that ``essentially distinct'' is a weaker condition than ``$\delta$-separated'' (at least for $c_0$ suitably large, which we will assume from now on).

In this article, we give an asymptotically sharp lower bound of $|\bigcup_{i\in I}T_i|$ for a collection of $N$ essentially distinct $\delta$-tubes in Section~\ref{bound}. We will see that the bound has a different expression for large $N$ and small $N$. Then we turn to the inverse problem: when the equality is approximately attained, what can we say about the configuration of these tubes? In Section~\ref{large}, we settle the situation when $N$ is large. In the rest of this article we focus on the more difficult situation when $N$ is small.

In Section~\ref{small}, we states our main rigidity result for small $N$. In Section~\ref{lemma}, we present some lemmas needed for its proof.

Let $\Omega=\Omega_n$ be the set of all lines in $\R^n$. Give $\Omega$ the canonical measure by identifying $\Omega$ with the product space $\S^{n-1}/\{\pm1\}\times\R^{n-1}$. The \textbf{X-ray transform} of a suitably defined function $f\colon\R^n\to\C$ is a function $Xf\colon\Omega\to\C$, $\ell\mapsto\int_\ell f$. In Section~\ref{X-ray} we reveal a close relationship between the $L^{n+1}$ norm of the $X$-ray transform of the indicator function of a set $E\subset\R^n$ and the convexity of $E$. This is an essential ingredient for our proof for the rigidity result.

On one hand, Christ \cite{ChristXray} proved (as a special case) that
\begin{equation}\label{eq:X}
||Xf||_{L^{n+1}(\Omega)}\le C||f||_{L^{(n+1)/2}(\R^n)}   
\end{equation}
for some constant $C>0$. On the other hand, it is a fact in integral geometry (see e.g. Ren \cite[(6.5.13)]{ren1994topics}) that $||X1_E||_{L^{n+1}(\Omega)}^{n+1}$ always takes the value $n(n+1)|E|^2/2$ when $E$ is convex, thus for $f=1_E$ the equality in \eqref{eq:X} is attained up to a constant. We will show that when $n\ge2$, in the converse direction, if we first suppose the equality in \eqref{eq:X} holds up to a constant for $f=1_E$, then $E$ is ``almost convex'' (see Theorem~\ref{thm:convex}). In particular, this gives raises to a reasonable measurement of convexity of sets (or more generally, functions) in $\R^n$ when $n\ge2$.

Finally, in Section~\ref{pf} we give the complete proof for the rigidity result for small $N$.

\vspace{15pt}

Throughout this article, $0<\delta<1/100$ is a small number, $n\ge2$ is an integer, $\T$ is a collection of essentially distinct $\delta$-tubes in $\R^n$, with $N=\#\T<\infty$.

We write $A\lesssim B$ for $A\le CB$ and $A\gtrsim B$ for $A\ge CB$ for some certain constant $C>0$ which may vary from line to line, and may depend on some prescribed constants (e.g. the dimension $n$) but not on any particular choice of variables (e.g. $\T$, $\delta$, $N$). We write $A\sim B$ for $A\lesssim B\lesssim A$.

For a set $E$, we write $1_E$ for its characteristic function. We use $B(x,r)=B_n(x,r)$ to denote the open ball in $\R^n$ with center $x$ and radius $r$. We use $e_1,\cdots,e_n$ to denote the coordinate vectors of $\R^n$. We use $|\cdot|$ to denote the Lebesgue measure in $\R^n$ or its induced surface/line measure for some surface/line in $\R^n$, which will be clear in the context. For a set $E\subset\R^n$ and a line $\ell\in\Omega_n$, we use $E_\ell$ to denote $E\cap\ell$.

A tube in $\R^n$ with radius $r$, height $h$, center $a$ and direction $e$ is the set defined by \eqref{tube} with $1/2$, $\delta/2$ replaced by $r$, $h/2$, respectively (so $T_e^\delta(a)$ is a tube with radius $\delta/2$, height $1$, center $a$ and direction $e$). For a direction $e$, the $\delta$-cap about $e$, denoted as $\c(e,\delta)$, is defined to be the set of all directions $e'\in\S^{n-1}/\{\pm1\}$ within angle $\delta$ of $e$.

\section{Sharp bound for the volume}\label{bound}
In this section, we give the asymptotically sharp lower bound for $|\bigcup_{T\in\T}T|$. We will use the ``bush'' argument first introduced by Bourgain \cite{bourga1991besicovitch}.
\begin{Thm}\label{thm:bound}
The following estimate is sharp up to constant factors:
\begin{equation}\label{eq:main}
    \left|\bigcup_{T\in\T}T\right|\gtrsim\begin{cases}
    \sqrt{N}\delta^{n-1}, &\text{ if }N\lesssim\delta^{2-2n}\\
    N\delta^{2n-2}, &\text{ if }N\gtrsim\delta^{2-2n}.
    \end{cases}
\end{equation}
\end{Thm}
\begin{pf}
We define a function $\mu$ on $\R^n$ by $\mu=\sum_{T\in\T}1_{T}$.

Denote $\nu=\max\mu$. Then we have the following estimate:
\begin{equation}\label{eq:e0}
    \left|\bigcup_{T\in\T}T\right|=\int_{\cup T}1\ge\int_{\cup T}\nu^{-1}\mu=\nu^{-1}\sum_{T\in\T}\int_{\R^n}1_{T}\sim \nu^{-1}N\delta^{n-1}.
\end{equation}

Consider any subcollection $\T'\subset\T$ whose tubes are within angle $\delta$ from each other. Then by an application of the pigeonhole principle, we see that every point $x\in\R^n$ is contained in $\lesssim1$ tubes in $\T'$. Since we can cover $\S^{n-1}/\{\pm1\}$ by $\sim\delta^{1-n}$ $\delta$-caps, we see that $\T$ can be written as the union of $\sim\delta^{1-n}$ such subcollection $\T'$. It follows that $\nu\lesssim\delta^{1-n}$.

Now \eqref{eq:e0} yields
\begin{equation}\label{eq:est1}
    \left|\bigcup_{T\in\T}T\right|\gtrsim N\delta^{2n-2}.
\end{equation}

On the other hand, choose a point $x\in\R^n$ with $\mu(x)=\nu$. Then there are $\nu$ tubes containing $x$. We can choose a suitably large constant $C\sim1$ so that for any two tubes $T_1,T_2$ containing $x$ with angle greater than $C\delta$, we always have $|T_1\cap T_2\cap\Sigma|=\emptyset$, where $\Sigma=\R^n\backslash B(x,1/4)$.

Since, by pigeonhole principle, there are $\lesssim1$ tubes in $\T$ containing $x$ with direction in any certain $C\delta$-cap, we can choose $\gtrsim\nu$ tubes in $\T$, say $T_1,\cdots,T_k$, with $C\delta$-separated directions. Then
\begin{equation}\label{eq:e3}
    \left|\bigcup_{T\in\T}T\right|\ge\left|\bigcup_{i=1}^k T_i\right|\ge\left|\bigcup_{i=1}^k(T_i\backslash\Sigma)\right|=\sum_{i=1}^k|T_i\backslash\Sigma|\gtrsim k\delta^{n-1}\gtrsim\nu\delta^{n-1}.
\end{equation}
Combining \eqref{eq:e0}\eqref{eq:e3}, we get
\begin{equation}\label{eq:est2}
    \left|\bigcup_{T\in\T}T\right|\gtrsim\sqrt{N}\delta^{n-1}.
\end{equation}

Now \eqref{eq:est1}\eqref{eq:est2} gives \eqref{eq:main}. We are left to show the estimate is sharp.

When $N\gtrsim\delta^{2-2n}$, we first put $\sim\delta^{1-n}$ disjoint $\delta$-tubes in the same direction together into a large tube with radius $\sim1$ and the same direction as the $\delta$-tubes. We fix a maximal collection of $\delta$-separated directions (which has cardinality $\sim\delta^{1-n}$). Then we put one copy of the large tube together with all $\delta$-tubes inside in each direction we have chosen, such that all the large tubes have the same center, say $O$. Now we have placed $\sim\delta^{2-2n}$ essentially distinct $\delta$-tubes, which can be covered by a ball of radius $\sim1$. We call this a standard configuration at $O$. Take $\sim N\delta^{2n-2}$ disjoint copy of some standard configurations, we have then constructed $\sim N$ $\delta$-tubes whose union has volume $\lesssim N\delta^{2n-2}$.

When $N\lesssim\delta^{2-2n}$, we first take $\sim\sqrt{N}$ disjoint $\delta$-tubes in the same direction together into a large tube with radius $\sim N^{1/(2n-2)}\delta$. We fix a maximal collection of $\delta$-separated directions restricted to a $N^{1/(2n-2)}\delta$-cap of $\S^{n-1}/\{\pm1\}$ (which has cardinality $\sim\sqrt{N}$). Then we put one copy of the large tube together with all $\delta$-tubes inside in each direction we have chosen, such that all the large tubes have the same center. Now we have placed $\sim N$ essentially distinct $\delta$-tubes, which can be covered by a tube of radius $\sim N^{1/(2n-2)}\delta$ (and height $2$, say). Thus the union of these $\delta$-tubes has volume $\lesssim(N^{1/(2n-2)}\delta)^{n-1}=\sqrt{N}\delta^{n-1}$.
\end{pf}

\section{Rigidity for large \texorpdfstring{$N$}{N}}\label{large}
Now we consider the inverse problem for Theorem~\ref{thm:bound}. That is, we suppose the equality holds (up to constant) in \eqref{eq:main}. We hope to show that $\T$ must look like the example we gave in Section~\ref{bound}. We shall first make this precise for $N\gtrsim\delta^{2-2n}$. We call the examples we constructed in the previous section "standard examples".

The equality holds means that there exists a constant $C>0$ independent of $\T$ such that
\begin{equation}\label{eq:large}
    \left|\bigcup_{T\in\T}T\right|\le CN\delta^{2n-2}.
\end{equation}

\begin{Def}
Let $T=T_{e}^{\delta}(a)$ be a $\delta$-tube, $O\in\R^n$ be any point. The \textbf{vertical distance} between $T$ and $O$ is $|(a-O)\cdot\hat e|$. The \textbf{horizontal distance} between $T$ and $O$ is $|(a-O)-((a-O)\cdot \hat e)\hat e|$. Here $\hat e$ is a preimage of $e$ in $\S^{n-1}$.
\end{Def}

\vspace{-.8em}

Thus, the $\delta$-tube $T=T_{e}^{\delta}(a)$ consists of all points that are within horizontal distance $\delta/2$ and vertical distance $1/2$ with respect to $T$.

\begin{Def}\label{d1}
For any $0<\epsilon_0<1$, $\lambda_0>0$, $O\in\R^n$, an \textbf{$(\epsilon_0,\lambda_0)$-good configuration at $O$} is a collection $\T_0$ of essentially distinct $\delta$-tubes which satisfies the following conditions:\vspace{-.8em}
\begin{enumerate}[(a)]
    \item All tubes in $\T_0$ are within horizontal distance $1/2$ and vertical distance $\epsilon_0$ with respect to $O$.
    \item There exists a collection $\{e_k\}_{k\in K}$ of $\delta$-separated directions with $|K|\ge\sqrt\lambda_0\delta^{1-n}$, such that $\T_0=\sqcup_{k\in K}\T_k$ for some $\{\T_k\}_{k\in K}$,
    where each $\T_k$ consists of at least $\sqrt\lambda_0\delta^{1-n}$ tubes, all of whose directions are in $\c(e_k,\epsilon_0\delta)$.
\end{enumerate}

We see that a $(\epsilon_0,\lambda_0)$-good configuration at $O$ looks very much like a subset of a certain standard configuration at $O$ we defined in the previous section. On one hand, $\epsilon_0$ measures how much a generic tube in a good configuration deviates from the one in the standard configuration: we can translate each tube vertically with respect to $O$ for no more than distance $\epsilon_0$ and then tilt each tube around its center for no more than angle $\epsilon_0\delta$ (and then possibly throw away a few tubes that intersect some other tubes in the same direction, but there will still be a positive amount of tubes remain) to get a subset of a standard configuration. On the other hand, $\lambda_0$ measures the portion of our good configuration possesses in the standard configuration: a $(\epsilon_0,\lambda_0)$-good configuration contains no less than $\lambda_0\delta^{2-2n}$ tubes, which has a portion $\sim\lambda_0$ of a standard configuration (this still holds even if we account for the tubes thrown away in the previous remark).

With this in mind, we see that the following theorem tells us that the standard examples we constructed in the previous section are essentially the only examples for equality to hold.

\begin{Thm}\label{thm:rigid_large}
There exists $\lambda>0$, such that for any $C,\epsilon_0>0$, there exists $\lambda_0>0$, such that for any $\T$ satisfying \eqref{eq:large}, we can extract at least $\lambda N\delta^{2n-2}$ $(\epsilon_0,\lambda_0)$-good configurations from $\T$, such that $(\bigcup_{T\in\T_1}T)\cap(\bigcup_{T\in\T_2}T)=\emptyset$ for any two different good configurations $\T_1,\T_2$ we extract.
\end{Thm}\vspace{-.8em}

Some remarks are in order.

\begin{Rmk}\hfill\vspace{-.8em}
\begin{enumerate}[(i)]
    \item We do not need to add the assumption $N\gtrsim\delta^{2-2n}$ since this is implied by \eqref{eq:main} and \eqref{eq:large}.
    \item By Definition~\ref{d1} we see that the total number of tubes we extract in the theorem is no less than $\lambda N\delta^{2n-2}\cdot \lambda_0\delta^{2-2n}=\lambda_0\lambda N$. This means a positive portion of tubes in $\T$ can be chosen to be an approximate copy of a subset of a (positive part of a) standard example.
    \item We cannot expect essentially better result than Theorem~\ref{thm:rigid_large}:\\
    For one thing, it is clear that we cannot require $\lambda_0$ (the portion of a good configuration in a standard configuration) to be independent of $C$ or $\epsilon_0$. For another thing, we can only guarantee to extract a positive portion of tubes but not an arbitrarily large portion of tubes to be contained in a standard example, even if we dropped the requirement that different good configurations are disjoint (which is nonessential). This is shown by the following example.
\end{enumerate}
\end{Rmk}

\begin{Ex}
Let $N_0\sim\delta^{2-2n}$. Use the example for small $N$ we constructed before, construct a configuration that contains $\sim N_0,N_0/2,N_0/2^2,\cdots,1$ $\delta$-tubes, respectively, such that the configurations are far away from each other. Call these configurations $\T_1,\cdots,\T_k$. Take $\T=\bigcup_{i=1}^k\T_i$. We see that $N\sim\delta^{2-2n}$, $|\bigcup_{T\in\T}T|\sim1$, so \eqref{eq:large} holds for some fixed $C>0$. 

To define rigidity, at least we have to fix a $\lambda_1>0$ which depends only on this fixed $C$, and say a collection of tubes is good if it (approximately) constitutes at least a portion $\lambda_1$ of some standard example. 

But after we fix this $\lambda_1$ and any positive integer $m$, for all sufficiently small $\delta$, we see that any subcollection $\T'\subset\T$ with $\#\T'\ge(1-2^{-m})N$ must intersect with $\gtrsim m$ of $\T_i$. If $\T'$ is a good collection of tubes, then there is a standard example $\T_0'$ such that $\T'$ constitutes a portion at least $\lambda_1$ of $\T_0'$. But by the construction of standard example we see that $\#\T_0'\gtrsim m\delta^{2-2n}$. So
\[
\delta^{2-2n}\sim N\ge(1-2^{-m})N\gtrsim\lambda_1m\delta^{2-2n}.
\]
In other words $m\lesssim\lambda_1^{-1}$. This means we cannot guarantee to extract an arbitrarily large portion of $\T$ to be (approximately) contained in a standard example.
\end{Ex}

We now prove the following two lemmas, from which Theorem~\ref{thm:rigid_large} follows.

\begin{Lem}
There exists $\lambda>0$ such that for any $C>0$, there exists $c>0$ such that for any $\T$ satisfying \eqref{eq:large}, there exist at least $\lambda N\delta^{2n-2}$ disjoint balls of radius $3$, each of which contains at least $c\delta^{2-2n}$ tubes in $\T$.
\end{Lem}
\begin{pf}
Let $\mu=\sum_{T\in\T}1_{T}$ denotes the multiplicity function of points in $\R^n$ with respect to $\T$. Then
\begin{equation}\label{eq:e_temp}
    \int_{\cup T}\mu=\sum_{T\in\T}\int_{\R^n}1_T\sim N\delta^{n-1}.
\end{equation}
On the other hand, we already see in the proof of Theorem~\ref{thm:bound} that $\mu\lesssim\delta^{1-n}$. This, together with \eqref{eq:large}\eqref{eq:e_temp}, implies that there exist some constants $\lambda_1>0$ independent of $C$ and $C_1=C_1(C)>0$ such that for
\[
M=\{\,x\in\bigcup_{T\in\T}T\colon\mu(x)\ge C_1\delta^{1-n}\,\},
\]
we have $|M|\ge2\lambda_1N\delta^{2n-2}$. 

Let
\[
A=\{\,x\in M\colon|B(x,1)\cap M|\le C_2\,\},\ B=M\backslash A,
\]
where $C_2>0$ is a constant to be chosen later.

Now we want to cover $A$ by a collection of balls $\{B(x_i,1)\}_{i=1}^k$ with $x_i\in A$ and $|x_i-x_j|\ge1$ for $i\ne j$.

To do this, select $x_i$ one by one. Having chosen $x_1,\cdots,x_{k'-1}$, we pick any $x_{k'}\in A\backslash(\bigcup_{i=1}^{k'-1}B(x_i,1))$ (if it exists). By pigeonhole principle, the condition $|x_i-x_j|\ge1$ implies every point in $\R^n$ is covered by $\lesssim1$ balls chosen. On the other hand, by the bush argument in the proof of Theorem~\ref{thm:bound} we see that 
\[
    \int_{\cup T}1_{B(x_i,1)}=|B(x_i,1)\cap(\bigcup_{T\in\T}T)|\gtrsim\mu(x_i)\delta^{n-1}\ge C_1,
\]
for all $1\le i\le k'$, thus
\[
k'C_1\lesssim\int_{\cup T}\sum_{i=1}^{k'}1_{B(x_i,1)}\lesssim\left|\bigcup_{T\in\T}T\right|\le CN\delta^{2n-2}.
\]
This means the selection process eventually stops at some
\[
    k\lesssim\frac{CN\delta^{2n-2}}{C_1}.
\]

With the chosen $\{B(x_i,1)\}_{i=1}^k$, we now have
\[
    |A|\le\sum_{i=1}^k|B(x_i,1)\cap M|\le kC_2\lesssim C_2\frac{CN\delta^{2n-2}}{C_1}.
\]
Thus we can choose some $C_2=C_2(C,C_1)=C_2(C)>0$ such that $|A|\le\lambda_1N\delta^{2n-2}$. Accordingly $|B|\ge\lambda_1N\delta^{2n-2}$.

Now find a collection of disjoint balls $\{B(x_j,3)\}_{j\in J}$ with $x_j\in B$. If $\{B(x_j,6)\}_{j\in J}$ does not cover $B$, then we can add more balls into this collection. Hence the collection can be chosen so that
\[
\#J\gtrsim|B|=\lambda_1N\delta^{2n-2}.
\]
Or written explicitly, $\#J\ge\lambda N\delta^{2n-2}$ for some $\lambda>0$ independent of $C$.

For each $j\in J$, by the choice of $M$ and $B$ we know that
\[
\int_{B(x_j,1)}\mu\ge\int_{B(x_j,1)\cap M}\mu\ge C_2C_1\delta^{1-n}.
\]
Since every $\delta$-tube intersects with $B(x_j,1)$ is contained in $B(x_j,3)$, the set $\T_j=\{\,T\in\T\colon T\subset B(x_j,3)\,\}$ satisfies
\[
\#\T_j\delta^{n-1}\sim\sum_{T\in\T_j}|T|\ge\int_{B(x_j,1)}\mu\ge C_1C_2\delta^{1-n}.
\]
Hence $\#\T_j\ge c\delta^{2-2n}$ for some $c=c(C)>0$. The proof is complete.
\end{pf}

\begin{Lem}\label{lem:center}
For any $c>0$, $0<\epsilon_0<1$, there exists $\lambda_0>0$ such that for any collection $\T_0$ of essentially distinct $\delta$-tubes contained in a ball $B$ of radius $3$, with cardinality at least $c\delta^{2-2n}$, there exists a subcollection of $\T_0$ which is a $(\epsilon_0,\lambda_0)$-good configuration.
\end{Lem}
\begin{pf}
For any $T\in\T_0$, let $T^\perp$ be the tube with radius $1/2$, height $2\epsilon_0$ which has the same center and direction as $T$.

Then $T^\perp\subset2B$, $\sum_{T\in\T}|T^\perp|\gtrsim\epsilon_0c\delta^{2-2n}$. Therefore there exists a point $O\in2B$ such that there are $\gtrsim\epsilon_0c\delta^{2-2n}$ tubes $T^\perp$ passing $O$.

Take $\T'=\{\,T\in\T_0\colon O\in T^\perp\,\}$. Then $\#\T'\gtrsim\epsilon_0c\delta^{2-2n}$ and that $(\T',O)$ satisfies (a) of Definition~\ref{d1}. Now it suffices to extract a subcollection of $\T'$ that satisfies (b).

Cover $\S^{n-1}/\{\pm1\}$ by a collection of $\delta$-caps $\{\Omega_i=\c(e_i,\delta)\}_{i\in I}$ such that each direction is covered by $\sim1$ caps in $\{\Omega_i\}_{i\in I}$. Let $\T_i=\{\,T\in\T'\colon T$ has direction in $\Omega_i\,\}$.

As remarked in the proof of Theorem~\ref{thm:bound}, each point in $\R^n$ is contained in $\lesssim1$ $\delta$-tubes in $\T'$ with direction in a certain $\delta$-cap. Now that all our tubes are contained in $B$, $|B|\sim1$, we thus know that $\#\T_i\lesssim\delta^{1-n}$. We also have $\sum_{i\in I}\#\T_i\sim\#\T'\gtrsim\epsilon_0c\delta^{2-2n}$, where $\#I\sim\delta^{1-n}$. From these we deduce that there exist $\lambda_1=\lambda_1(c,\epsilon_0)$, $\lambda_2=\lambda_2(c,\epsilon_0)>0$ such that the set
\[
I'=\{\,i\in I\colon\#\T_i\ge\lambda_2\delta^{1-n}\,\}
\]
has cardinality at least $\lambda_1\delta^{1-n}$.
Decrease $\lambda_1$ if necessary, we may also assume $\{\,e_i\colon i\in I'\,\}$ has pairwise angle larger than $4\delta$.

For each $i\in I'$, we further cover $\Omega_i$ by $\sim\epsilon_0^{1-n}$ $\epsilon_0\delta/2$-caps. By pigeonhole principle we can find $\gtrsim\epsilon_0^{n-1}\lambda_1\delta^{1-n}$ (written more explicitly, $\ge c'\epsilon_0^{n-1}\lambda_1\delta^{1-n}$) tubes in $\T_i$ with direction in one of those $\epsilon_0\delta/2$-caps. Take $\lambda_0=\min\{(c'\epsilon_0^{n-1}\lambda_1)^2,\lambda_2^2\}$ and the proof is complete.
\end{pf}

\end{Def}

\section{Rigidity for small \texorpdfstring{$N$}{N}}\label{small}
Now we turn to the situation when $N$ is small. In this case, the equality of \eqref{eq:main} takes the form
\begin{equation}\label{eq:small}
    \left|\bigcup_{T\in\T}T\right|\le C\sqrt{N}\delta^{n-1}
\end{equation}
for some fixed constant $C>0$.

We hope to prove similar result as in the previous section, i.e. our example constructed earlier is unique in some sense. However, it is not true that we can extract a certain amount of tubes in $\T$ that resembles a subset of the ``standard example''. In fact, examples with different outlook do exist.

\begin{Ex}
When $N\lesssim\delta^{-2d}$ for some $2\le d\le n$, we first employ the ``standard example'' we have constructed when the dimension of the whole space is $d$ to obtain a sharp example of a collection $\T'$ of essentially distinct $\delta$-tubes in $\R^d$. For each $T\in\T'$, take $\widetilde T\subset\R^n$ to be the $\delta$-tube whose intersection with $\R^d\times\{0\}^{n-d}\simeq\R^d$ is $T$, and whose center and direction are both contained in $\R^d\times\{0\}^{n-d}$. Then $\T=\{\,\widetilde T\colon T\in\T'\,\}$ gives a sharp example for dimension $n$.
\end{Ex}

We now have $n-2$ seemingly different examples which satisfy \eqref{eq:small} for some large $C$ (at least when $N\lesssim\delta^{-4}$). There is a better way to describe these examples.

\begin{Ex}
Suppose $N\lesssim\delta^{-2d}$, so that $N^{1/2d}\delta\lesssim1$.

Let $E=[0,N^{1/2d}\delta]^{d}\times[0,\delta]^{n-1-d}\subset\R^{n-1}$. Put as many essentially distinct $\delta$-tubes into $E\times[0,2]$ as possible.

We see that there are $\sim\sqrt{N}$ $\delta$-separated directions in each of which we can put $\sim\sqrt{N}$ essentially distinct tubes. Now the total number of $\delta$-tubes we put is $\sim N$ and their union has volume $\lesssim\sqrt{N}\delta^{n-1}$.
\end{Ex}

Based on this description, we can make two guesses. First, all examples satisfying \eqref{eq:small} arise in similar way. i.e. there is a ``good'' set $E\subset\R^{n-1}$ with $|E|\sim\sqrt{N}\delta^{n-1}$ such that most of the tubes are contained in a set congruent to $E\times[0,2]$. Second, every example arises this way does contain $\sim N$ tubes, thus is a sharp example. We shall now make these precise.

What does ``good'' mean? Since our tubes has radius $\delta/2$, it is natural to require $E$ to be \textit{(a) $k\delta$-discretized} (i.e. is a union of balls of radius $k\delta$) for some constant $k>0$. But this is obviously not enough. For example, we take $E$ to be an annulus with width $k\delta$ and volume $\sqrt{N}\delta^{n-1}$, then for $n\ge3$, it is not possible to put as many as $\sim N$ essentially distinct $\delta$-tubes inside $E\times[0,2]$. A natural guess is that $E$ should also be \textit{(b) convex}. It turns out that these two conditions, plus the assumption that $E$ \textit{(c) does not have too large diameter}, are sufficient for both the two guesses.

We first deal with the second guess. That is to say, \textit{$(a)(b)(c)$} implies ``good''.

\begin{Thm}\label{thm:convex_is_good}
Let $E\subset\R^{n-1}$ be a convex $9\delta$-discretized set with $|E|\sim\sqrt{N}\delta^{n-1}$ and $\diam(E)\le1$. Then we can put $\sim N$ essentially distinct $\delta$-tubes into $E\times[0,2]$.
\end{Thm}

\begin{Rmk}
Of course, our proof presented below is not the best way to solve this problem (see Remark~\ref{rmk:convex}). However, we do want to show some continuization argument which turns discrete problems into continuous ones. Similar technique will be needed in the proof of Theorem~\ref{thm:small}, and we shall not repeat it again then.
\end{Rmk}

\begin{pf}
Take a collection of $\delta$-separated directions $\I=\{e_i\}$ such that $\c(e_i,\delta/2)$ are pairwise disjoint but $\bigcup_{e_i\in\I}\c(e_i,\delta)=\S^{n-1}/\{\pm1\}$. Then $\#\I\sim\delta^{1-n}$.

For any direction $e$, let $\#'e$ denote the maximal number of essentially distinct $\delta$-tubes with direction $e$ that can be placed into $E\times[0,2]$, $\#e$ denote the maximal number of essentially distinct tubes with radius $3\delta/2$, height $3/2$ and direction $e$ that can be placed into $E\times[0,2]$.

For any $e$, we can find $e_i\in\I\cap\c(e,\delta)$. Since any $\delta$-tube with direction $e_i$ is contained in a tubes with radius $3\delta/2$, height $3/2$ and direction $e$, we see that $\#'e_i\ge\#e$. Thus
\begin{equation}\label{eq:approx_1}
    \sum_{e_i\in\I}\#'e_i\gtrsim\delta^{1-n}\sum_{e_i\in\I}\int_{\c(e_i,\delta)}\#e\ge\delta^{1-n}\int_{\S^{n-1}/\{\pm1\}}\#e=\delta^{1-n}\int_{(\S^{n-1})_+}\#\hat e,
\end{equation}
where $(\S^{n-1})_+=\S^{n-1}\cap(\R^{n-1}\times(0,\infty))$, $\hat e$ is the preimage of $e$ in $(\S^{n-1})_+$, $\#\hat e=\#e$.

Now we replace $E$ by $\frac{2}{3}E$, and we redefine $\#\hat e$ to be the maximal number of essentially distinct $2\delta$-tubes with direction $e$ that can be placed into $E\times[0,4/3]$. Then $\#\hat e$ agrees with the previous definition. It suffices to show that
\begin{equation}\label{eq:ets}
    \int_{(\S^{n-1})_+}\#\hat e\gtrsim N\delta^{n-1}.
\end{equation}
We may further replace the $2\delta$-tubes in the definition of $\#\hat e$ by $\delta$-tubes, and by doing a rescaling, we are left to show \eqref{eq:ets} provided that $E\subset\R^{n-1}$ is a convex $3\delta$-discretized set with $|E|\sim\sqrt{N}\delta^{n-1}$ and $\diam(E)\le1$.

Let $E'=\{\,x\in E\colon d(x,\partial E)>2\delta\,\}$. Then $|E'|\sim|E|$, $|E'|$ is convex and $\delta$-discretized.

Write $\hat e=\sin\theta\cdot\xi+\cos\theta\cdot e_n$ where $\xi\in\R^{n-1}\times\{0\}$ is a unit vector, $\theta\in[0,\pi/2)$, and $e_n$ is the $n$-th coordinate vector. Denote $\#\hat e$ by $\#\xi_\theta$. Then
\begin{equation}\label{eq:1}
    \int_{(\S^{n-1})_+}\#\hat e\sim\int_{\S^{n-2}}\int_0^{\frac{\pi}{2}}\theta^{n-2}\#\xi_\theta
\end{equation}

For fixed $\xi$, divide $\R^n$ into $2$-slices $\R^n=\bigcup_{t}S_t$ whose interiors are pairwise disjoint, where each $2$-slice $S_t$ is congruent to $\R^{2}\times[0,\delta]^{n-2}$, such that the two dimensions correspond to $\R^2$ are spanned by $\xi$ and $e_n$. Denote $A_t=S_t\cap(E\times\{0\})$, $A'_t=S_t\cap(E'\times\{0\})$.

Let $\Pi:\R^n\to\R^{n-1}\times\{0\}\simeq\R^{n-1}$ be the projection map. For a $\delta$-tube $T$ with direction $e$, $\Pi(T)$ is contained in a tube of dimension $n-1$ with radius $\delta/2$ and height $\sin\theta+\delta\le\theta+\delta$. Thus a translation of it can be contained in $A_t$ if $|(A_t')_{\ell_t}|(=|(E')_{\ell_t}|)\ge\theta$, where $\ell_t$ is the central line in the $1$-slice $P_t=\Pi(S_t)$. Moreover, in the vertical direction (the last coordinate), $T$ has width $\lesssim\delta/\theta$ above each point in $\R^{n-1}\times\{0\}$ and a total width in the vertical direction of at most $1+\delta<101/100$. Thus if a translation of $\Pi(T)$ can be contained in $A_t$, then a total of $\gtrsim\theta/\delta$ essentially distinct $\delta$-tubes with direction $e$ can be placed into $S_t\cap(E\times[0,4/3])$. The above argument shows that
\[
\#\xi_\theta\gtrsim\sum_{t\colon|(E')_{\ell_t}|\ge\theta}\frac{\theta}{\delta}.
\]
(By $\diam E\le1$ one can check that this is actually an equality rather than inequality. The other direction is used for the proof of Theorem~\ref{thm:small}.)

Thus
\begin{equation}\label{eq:2}
    \int_0^{\frac{\pi}{2}}\theta^{n-2}\#\xi_\theta d\theta\gtrsim\sum_t\int_0^{|(E')_{\ell_t}|}\frac{\theta^{n-1}}{\delta}\sim\frac{1}{\delta}\sum_t|(E')_{\ell_t}|^n.
\end{equation}
Here we used the fact that $|(E')_{\ell_t}|\le\diam(E')\le\diam(E)<\pi/2$.

Now let $E''=\{\,x\in E'\colon d(x,\partial E')>\delta/2\,\}$. Then $E''$ is convex and $|E''|\sim|E|\sim\sqrt{N}\delta^{n-1}$. We have
\begin{equation}\label{eq:3}
    |(E')_{\ell_t}|\ge|(E'')_\ell|\text{ for all }\ell\subset P_t.
\end{equation}
Now \eqref{eq:1}\eqref{eq:2}\eqref{eq:3} yield
\begin{equation}\label{eq:approx_2}
    \int_{(\S^{n-1})_+}\#\hat e\gtrsim\frac{1}{\delta}\int_{\S^{n-2}}\sum_{t}\delta^{2-n}\int_{\ell\subset P_t}|(E'')_{\ell}|^n=\delta^{1-n}\int_{\Omega_{n-1}}|(E'')_\ell|^n,
\end{equation}
where $\Omega_{n-1}$ denotes the set of all lines in $\R^{n-1}$ with the natural measure.

But it is a known fact (see Ren \cite[(6.5.13)]{ren1994topics}) that for any convex set $K\subset\R^{m}$, we have
\[
\int_{\Omega_m}|K_\ell|^{m+1}=\frac{m(m+1)}{2}|K|^2.
\]
Thus \eqref{eq:ets} follows from \eqref{eq:approx_2}, and the proof is complete.
\end{pf}

Proving the other guess we made is much harder. We have only an (approximate) equality on our hand but are asking for rather rigid structural property for a given sharp example. We first state out this rigidity result, which is one of the main theorems in our paper.

\begin{Thm}\label{thm:small}
Suppose $\T$ satisfies \eqref{eq:small}, then we can find a convex $9\delta$-discretized set $E\subset\R^{n-1}$ with $|E|\sim\sqrt{N}\delta^{n-1}$ and $\diam(E)\le1$, such that there are $\sim N$ tubes in $\T$ contained in a set that is congruent to $E\times[0,2]$. Here the implicit constants depend only on $C,n$.
\end{Thm}

\begin{Rmk}\label{rmk:convex}
There is an alternative way to define ``good''. We may define a ``good'' set to be any box (resp. ellipsoid) in $\R^n$ with all side length (resp. axis length) belonging to $[k\delta,1]$ for a suitable constant $k$. The equivalence between this definition and the one we stated follows from the fact that for any convex set $E\subset\R^n$ with non-empty interior, there exists homothetic boxes (resp. ellipsoids) $B_-,B_+$ with $B_-\subset E\subset B_+$ and $|B_-|\sim|B_+|\sim|E|$ (see \cite{GRUBER1993319}). Nevertheless, we stick with our previous definition.
\end{Rmk}

We extract several technical tools for our proof of Theorem~\ref{thm:small} and present them in the next two sections. Then we give the complete proof in Section~\ref{pf}.

\section{Lemmas}\label{lemma}

\begin{Lem}\label{lem:large_simplex}
For every $c>0$ and positive integer $n$, there exists $\lambda,c'>0$ such that the following holds:\\
For any measure space $(E,\M,\mu)$ with $0<\mu(E)<\infty$, any measurable set $S\subset\R^n$ with $|S|>0$, and an assignment $x\mapsto E_x\in\M$ with $\mu(E_x)\ge c\mu(E)$ for each $x\in S$, there exists $x_0,\cdots,x_n\in S$ such that:\vspace{-.8em}
\begin{enumerate}[(a)]
    \item $\mu(\bigcap_{i=0}^n E_{x_i})\ge c'\mu(E)$;
    \item $\vol(x_0,\cdots,x_n)\ge\lambda|S|$.
\end{enumerate}\vspace{-.8em}
Where $\vol(x_0,\cdots,x_n)$ denotes the volume of the $n$-simplex with vertices $x_0,\cdots,x_n$.
\end{Lem}
\begin{pf}
For $m=0,1,\cdots,n-1$, we denote the $m$-hyperplane determined by $m+1$ points $x_0,x_1,\cdots,x_m\in\R^n$ in general position by $\alpha_m(x_0,\cdots,x_m)$. For any $m$-hyperplane $\alpha_m$ and $r>0$, we define $D(\alpha_m,r)\subset\R^n$ to be the $r$-neighborhood of $\alpha_m$. Under these notations, the ball $B(x,r)$ in $\R^n$ can also be written as $D(x,r)=D(\alpha_0(x),r)$.

Below we give an algorithm to explicitly find the desired points $x_0,\cdots,x_n$.

Set $I=(1,\cdots,1)$ (a total $n$ of $1$'s), $E^I=E$, $S^I=S$, $E_x^I=E_x$, $c^I=c$. We begin the following process:

{\addtolength{\leftskip}{5mm}
\textit{Step $I=(i_1,\cdots,i_n)$}:\\
Pick any $x_0\in E^I$. Set $k=0$. We begin the following subprocess:

{\addtolength{\leftskip}{5mm}
\textit{Substep $k$}:\\
If $k=n$ then stop the subprocess. Otherwise, choose $r_{k+1}>0$ such that $D_{k+1}=D(\alpha_k(x_0,\cdots,x_k),r_{k+1})$ satisfies
\[
|S^I\cap(\bigcap_{j=1}^{k+1}D_j)|=\frac{|S^I|}{2^{k+1}}.
\]
If there exists $x_{k+1}\in S^I\cap(\bigcap_{j=1}^kD_j)\backslash D_{k+1}$ such that
\begin{equation}\label{eq:cap_E}
    \mu(\bigcap_{j=0}^{k+1}E_{x_j}^I)\ge\frac{c^I\mu(E)}{((\cdots((1+i_1)^2+i_2)^2+\cdots+i_k)^2+i_{k+1})^2},
\end{equation}
then goto \textit{Substep $k+1$}, otherwise stop the subprocess.

}
\vspace{10pt}
The subprocess must stop at some $k\le n$. If $k=n$, then we have selected $x_0,\cdots,x_n$, and we stop the whole process. Otherwise, for all $x\in S^I\cap(\bigcap_{j=1}^kD_j)\backslash D_{k+1}$ we have
\[
\mu((\bigcap_{j=0}^kE_{x_j}^I)\cap E_x^I)<\frac{c^I\mu(E)}{((\cdots((1+i_1)^2+i_2)^2+\cdots+i_k)^2+i_{k+1})^2}.
\]

Now we define
\[
I'=(i_1,\cdots,i_k,i_{k+1}+1,1,\cdots,1),
\]
\[
E^{I'}=E^I\backslash(\bigcap_{j=0}^kE_{x_j}^I),
\]
\[S^{I'}=S^I\cap(\bigcap_{j=1}^kD_j)\backslash D_{k+1},
\]
\[E_x^{I'}=E_x^I\backslash(\bigcap_{j=0}^kE_{x_j}^I),\text{ for all }x\in S^{I'},
\]
\begin{equation}\label{eq:c_t}
   c^{I'}=\left(1-\frac{1}{((\cdots((1+i_1)^2+i_2)^2+\cdots+i_k)^2+i_{k+1})^2}\right)c^I. 
\end{equation}
And we goto \textit{Step $I'$}.

}
\vspace{10pt}
Now let us analyze the above process. By an induction we see that for each step $I$, we have
\[
E_x^I\subset E^I,\ \mu(E_x^I)\ge c^I\mu(E)\text{ for all }x\in S^I.
\]
Every time when we move from \textit{Step $I=(i_1,\cdots,i_n)$} to \textit{Step $I'=(i_1,\cdots,i_k,i_{k+1}+1,1,\cdots,1)$}, by choice of $x_k$ we have
\begin{equation}\label{eq:E_t}
    \mu(E^{I'})\le\mu(E^I)-\mu(\bigcap_{j=0}^kE_{x_j}^I)\le\mu(E^I)-\frac{c^{I}\mu(E)}{(\cdots((1+i_1)^2+i_2)^2+\cdots+i_k)^2}.
\end{equation}
(The denominator in the last term is $1$ if $k=0$.)

Also,
\begin{equation}\label{eq:S_t}
|S^{I'}|=|S^I\cap(\bigcap_{j=1}^kD_j)|-|S^I\cap(\bigcap_{j=1}^{k+1}D_j)|=\frac{|S^I|}{2^{k+1}}\ge\frac{|S^I|}{2^n}.
\end{equation}
From \eqref{eq:c_t} we know that for any step I,
\begin{align}\label{eq:c^I}
c^I\ge&\ c\prod_{i_1=1}^\infty\left(\left(1-\frac{1}{(i_1+1)^2}\right)\prod_{i_2=1}^\infty\left(\left(1-\frac{1}{((1+i_1)^2+i_2)^2}\right)\cdots\vphantom{\prod_{i_n=1}^\infty\left(1-\frac{1}{((\cdots((1+i_1)^2+i_2)^2+\cdots+i_n)^2+i_{n+1})^2}\right)}\right.\right.\nonumber\\ &\left.\left.\prod_{i_n=1}^\infty\left(1-\frac{1}{((\cdots((1+i_1)^2+i_2)^2+\cdots+i_{n-1})^2+i_n)^2}\right)\cdots\right)\right)\nonumber\\
=&\ c\prod_{k=1}^n\left(\prod_{i_1=1}^\infty\cdots\prod_{i_k=1}^\infty\left(1-\frac{1}{((\cdots((1+i_1)^2+i_2)^2+\cdots+i_{k-1})^2+i_k)^2}\right)\right).
\end{align}
Using the simple fact $\prod_{i=1}^\infty(1-\frac{1}{(M+i)^2})=\frac{M}{M+1}>1-\frac{1}{M}$, we know that
\begin{align*}
    &\prod_{i_1=1}^\infty\cdots\prod_{i_k=1}^\infty\left(1-\frac{1}{((\cdots((1+i_1)^2+i_2)^2+\cdots+i_{k-1})^2+i_k)^2}\right)\\
    >&\prod_{i_1=1}^\infty\cdots\prod_{i_{k-1}=1}^\infty\left(1-\frac{1}{(\cdots((1+i_1)^2+i_2)^2+\cdots+i_{k-1})^2}\right)\\
    >&\cdots\\
    >&\prod_{i_1=1}^\infty(1-\frac{1}{(i_1+1)^2})\\
    =&\ \frac{1}{2}.
\end{align*}
Thus \eqref{eq:c^I} gives
\begin{equation}\label{eq:c}
c^I\ge\frac{c}{2^n}.    
\end{equation}

Next, we show that the process would eventually stop, and the total number of steps is bounded above by some constant $K=K(c,n)>0$.

We use induction on $k$ to show that the quantity $i_k$ in any step $I=(i_1,\cdots,i_n)$ is bounded above by a constant $K_k=K_k(c,n)$.

For $0\le k\le n-1$, suppose the the bound $K_{k_0}$ exists for all $k_0\le k$, consider $k+1$. Denote $I_1=(i_1,\cdots,i_k,1,\cdots,1)$, $I_2=(i_1,\cdots,i_{k+1},1,\cdots,1)$, \eqref{eq:E_t}\eqref{eq:c} give that
\begin{align*}
    0\le&\ \mu(E^{I_2})\\
    \le&\ \mu(E^{I_1})-\sum_{j=1}^{i_{k+1}-1}\frac{1}{(\cdots((1+i_1)^2+i_2)^2+\cdots+i_k)^2}\frac{c}{2^n}\mu(E)\\
    \le&\ \mu(E)-(i_{k+1}-1)\frac{1}{(\cdots((1+K_1)^2+K_2)^2+\cdots+K_k)^2}\frac{c}{2^n}\mu(E).
\end{align*}
Thus it suffices to take
\[
K_{k+1}=1+\frac{2^n(\cdots((1+K_1)^2+K_2)^2+\cdots+K_k)^2}{c}.
\]
(The numerator in the last term is $2^n$ if $k=0$.)

The induction is complete, and the assertion of the boundedness of the total number of steps follows.

Now suppose the whole process stop at \textit{Step $I_0=(i_1,\cdots,i_n)$}. We shall show that the $n+1$ points $x_0,\cdots,x_n$ satisfy our desired properties \textit{(a)(b)} for some constants $c',\lambda>0$.

\textit{(a)}: By \eqref{eq:cap_E}\eqref{eq:c},
\begin{align*}
    \mu(\bigcap_{i=0}^n E_{x_i})\ge&\ \frac{1}{(\cdots((1+i_1)^2+i_2)^2+\cdots+i_n)^2}\frac{c}{2^n}\mu(E)\\
    \ge&\ \frac{1}{(\cdots((1+K_1)^2+K_2)^2+\cdots+K_n)^2}\frac{c}{2^n}\mu(E).
\end{align*}

\textit{(b)}: By \eqref{eq:S_t} and the boundedness of total number of steps,
\[
|S^I|\ge\frac{|S|}{2^{(K-1)n}}.
\]
By the construction of $r_k$, $D_k$, $x_k$, we see that
\begin{equation}\label{eq:cap_D}
    |\bigcap_{j=1}^nD_j|\ge|S^I\cap(\bigcap_{j=1}^nD_j)|=\frac{|S^I|}{2^n}\ge\frac{|S|}{2^{Kn}},
\end{equation}
and that for each $k=0,\cdots,n-1$, the distance $d_{k+1}$ from $x_{k+1}$ to the $k$-hyperplane $\alpha_k(x_0,\cdots,x_k)$ satisfies
\begin{equation}\label{eq:d}
    d_{k+1}\ge r_{k+1}.
\end{equation}
But $\bigcap_{j=1}^nD_j$ is contained in a box of size $2r_1\times\cdots\times2r_n$, thus \eqref{eq:cap_D} further yields
\[
\prod_{k=1}^nr_k\ge\frac{|S|}{2^{(K+1)n}},
\]
and then \eqref{eq:d} yields
\[
\vol(x_0,\cdots,x_n)=\frac{1}{n!}\prod_{k=1}^nd_k\ge\frac{|S|}{2^{(K+1)n}n!}.
\]
The proof is complete.
\end{pf}

\begin{Lem}\label{lem:small_int}
Let $E\subset\R^n$ be a measurable set with $d=\diam(E)$. For $t\in\R$, denote $E_t=E\cap(\bigcap_{i=1}^n(E+te_i))$. Then
\[
\int_\R|E_t|dt\le2(d^{n-1}|E|^{2n})^{\frac{1}{2n-1}}.
\]
\end{Lem}
\begin{pf}
Define 
\[
\phi\colon\R^n\to\R,\ x\mapsto\int_\R\prod_{i=1}^n1_E(x+te_i)dt.
\]
Then $\phi(x)\le\int_\R1_E(x+te_1)\le\diam(E)\le d$ for all $x\in\R^n$. Denote $E'=\{\,x\in E\colon\phi(x)\ge c\,\}$, where $c>0$ is a constant to be chosen later.

Let $\Pi_i\colon\R^n\to e_i^\perp=\R^{i-1}\times\{0\}\times\R^{n-i}$ be the projection map along $e_i$, $i=1,\cdots,n$. For any $\tilde x\in\Pi_i(E')$, take its preimage $x_0\in\R^n$, we have
\[
\int_\R1_E(\tilde x+te_i)dt=\int_\R1_E(x_0+te_i)dt\ge\phi(x_0)\ge c.
\]
Thus
\[
|E|=\int_{e_i^\perp}\int_\R1_E(\tilde x+te_i)dtd\tilde x\ge\int_{\Pi_i(E')}\int_\R1_E(\tilde x+te_i)dtd\tilde x\ge c|\Pi_i(E')|.
\]
Apply the isoperimetric inequality to $E'$, we get
\[
|E|^n\ge c^n\prod_{i=1}^n|\Pi_i(E')|\ge c^n|E'|^{n-1}.
\]
In other words
\[
|E'|\le\left(\frac{|E|}{c}\right)^{\frac{n}{n-1}}.
\]
Now we have estimate
\begin{align*}
    \int_\R|E_t|dt=&\int_\R\int_E\prod_{i=1}^n1_E(x+te_i)dxdt\\
    =&\int_E\phi(x)dx\\
    =&\int_{E'}\phi(x)dx+\int_{E\backslash E'}\phi(x)dx\\
    \le&\ |E'|d+|E|c\\
    \le&\ \frac{d|E|^\frac{n}{n-1}}{c^\frac{n}{n-1}}+c|E|\\
    =&\ 2(d^{n-1}|E|^{2n})^{\frac{1}{2n-1}},
\end{align*}
where we took $c=(d^{n-1}|E|)^{1/(2n-1)}$.
\end{pf}

\begin{Lem}\label{lem:arithmetic}
Let $(Z,+)$ be a free abelian group, $A,B\subset Z$, $G\subset A\times B$, $C=\{a+b\colon(a,b)\in G\}$, such that $\#A,\#B,\#C\le N_0$, and that the map $G\to Z$, $(a,b)\mapsto a-b$ maps at most $M$ points to one point. Then $\#G\le M^{1/6}N_0^{11/6}$.
\end{Lem}

Katz, Tao proved this result for $M=1$. But there is no significant difference in the proof for this general version, as remarked by Oberlin \cite{Oberlin2010}. We shall avoid unnecessary repetition here and refer the readers to \cite{katz1999bounds}.

\section{X-ray transform and convexity}\label{X-ray}

We propose a new measurement for the convexity of sets in $\R^n$, which follows from a more general definition given for suitably defined functions in $\R^n$. We shall see that our definitions (especially the one for sets) have many nice properties. For some other measurements of convexity, one may consult \cite{MANILEVITSKA199319}. In particular, our definition is close to that of Beer \cite{beer1974index}, but has its own advantage.

In this section, we will always assume $n\ge2$.

\begin{Def}
The \textbf{convexity index} of a function $f\colon\R^n\to\C$ is defined to be
\[
c(f)=c_n(f):=\frac{2}{n(n+1)}\left(\frac{||Xf||_{L^{n+1}(\Omega_n)}}{||f||_{L^{(n+1)/2}(\R^n)}}\right)^{n+1}.
\]
The \textbf{convexity index} of a set $E\subset\R^n$ is defined to be
\[
c(E)=c_n(E):=c(1_E)=\frac{2}{n(n+1)|E|^2}\int_{\Omega_n}|E_\ell|^{n+1}d\ell,
\]
where $\Omega_n$ is the set of all lines in $\R^n$ with the natural measure.
\end{Def}

\begin{Rmk}
Although our definition of convexity index has several advantages over the one defined by Beer \cite{beer1974index}, it fails to measure convexity when $n=1$.
\end{Rmk}

\begin{Prop}\label{prop:index}
The convexity index of a function $f\colon\R^n\to\C$ (resp. a set $E\subset\R^n$) is invariant under non-singular affine transformations or rearrangement on a set of measure zero. Also, $c(f)$ is invariant under multiplication by nonzero constants.
\end{Prop}
\begin{Rmk}
The Beer convexity index (defined only for sets) is also invariant under non-singular affine transformations, but it is sensitive to changes on a set of measure zero. From this point of view, our definition of convexity index is more natural under the analytical sense.
\end{Rmk}
\begin{pf}
Let $A$ be any non-singular affine transformation in $\R^n$. Write $f_A(x)=f(Ax)$. Then we have
\begin{equation}\label{eq:f_A}
||f_A||_{L^{(n+1)/2}}^{(n+1)/2}=\int_{\R^n} f(Ax)^{(n+1)/2}dx=|\det A|^{-1}||f||_{L^{(n+1)/2}}^{(n+1)/2}.
\end{equation}
For any line $\ell\in\Omega$ with direction $\xi$, we have
\[
Xf_A(\ell)=\int_\ell f(Ax)dx=\int_{A\ell} f(y)d(A^{-1}y)=|A\xi|^{-1}Xf(A\ell).
\]
For $\xi\in\S^{n-1}/\{\pm1\}$, $x\in\R^n$, let $\xi^\perp$ denotes an orthogonal complement of $\xi$ in $\R^n$, $\ell(\xi,x)$ denotes the line passing $x$ with direction $\xi$. Then
\begin{align}\label{eq:Xf_A}
    ||Xf_A||_{L^{n+1}}^{n+1}=&\int_{\S^{n-1}/\{\pm1\}}\int_{\xi^\perp}Xf_A(\ell(\xi,x))^{n+1}dxd\xi\nonumber\\
    =&\int_{\S^{n-1}/\{\pm1\}}|A\xi|^{-(n+1)}\int_{\xi^\perp}Xf(\ell(\frac{A\xi}{|A\xi|},Ax))^{n+1}dxd\xi\nonumber\\
    =&\int_{\S^{n-1}/\{\pm1\}}|A\xi|^{-n}|\det A|^{-1}\int_{(\frac{A\xi}{|A\xi|})^\perp}Xf(\ell(\frac{A\xi}{|A\xi|},y))^{n+1}dyd\xi\nonumber\\
    =&|\det A|^{-2}\int_{\S^{n-1}/\{\pm1\}}\int_{\zeta^\perp}Xf(\ell(\zeta,y))^{n+1}dyd\zeta\nonumber\\
    =&|\det A|^{-2}||Xf||_{L^{n+1}}^{n+1}.
\end{align}
Now \eqref{eq:f_A}\eqref{eq:Xf_A} show $c(f_A)=c(f)$.

Suppose $f_0=f+g$ where $g=0$ almost everywhere. Then for fixed $\xi\in\S^{n-1}/\{\pm1\}$, $\int_{\ell(\xi,x)} g$ exists and equals to zero for almost every $x\in\xi^\perp$, so
\[
\int_{\xi^\perp}Xf_0(\xi,x)^{n+1}dx=\int_{\xi^\perp}Xf(\xi,x)^{n+1}dx.
\]
Thus $||Xf_0||_{L^{n+1}}=||Xf||_{L^{n+1}}$. On the other hand, it is clear that $||f_0||_{L^{(n+1)/2}}=||f||_{L^{(n+1)/2}}$. It follows that $c(f_0)=c(f)$.

The statements for sets follow from those for functions.

The invariance of $c(f)$ under multiplication by nonzero constant is obvious.
\end{pf}

\begin{Thm}\label{thm:index}\hfill\vspace{-.8em}
\begin{enumerate}[(1)]
    \item For any $f\colon\R^n\to\C$ whose convexity index is defined, we have $0\le c(f)\lesssim1$, where the implicit constant depends only on $n$.
    \item For any set $E\subset\R^n$ whose convexity index is defined, we have $0\le c(E)\le1$, and the equality on the right holds if and only if $E$ is convex up to rearrangement of points in a set of measure zero.
\end{enumerate}
 
\end{Thm}\vspace{-.8em}
The statement for functions follows from the following $L^p$ estimate for the X-ray transform:
\[
||Xf||_{L^{n+1}(\Omega)}\lesssim||f||_{L^{(n+1)/2}(\R^n)},
\]
which was first proved (as a special case) by Christ \cite{ChristXray}.

As for the statement for sets, we shall follow the computation of the equality case for convex sets (Ren \cite[(6.5.13)]{ren1994topics}) and give a proof for our more general statements.

\begin{pf}
Using a change of variable formula (see Ren \cite[(6.2.37)]{ren1994topics}), we have
\begin{equation}\label{eq:E^2}
    |E|^2=\int_{\R^n}\int_{\R^n}1_E(x)1_E(y)dxdy=\int_\Omega\int_\ell\int_\ell|t_1-t_2|^{n-1}1_E(t_1)1_E(t_2)dt_1dt_2d\ell.
\end{equation}
For a fixed $\ell$, let us identify $\ell$ with $\R$. For $t\in\ell$, define $L(t)=|(-\infty,t]\cap E_\ell|$. Then $L\colon\ell\to[0,|E_\ell|]$ is monotonically increasing and we can choose a monotonically increasing $I\colon(0,|E_\ell|)\to\ell$ such that $L\circ I=\id_{(0,|E_\ell|)}$. Then
\begin{equation}\label{eq:Lip}
    |I(t_1)-I(t_2)|\ge|t_1-t_2|\text{ for all }t_1,t_2\in(0,|E_\ell|),
\end{equation}
and thus
\begin{align}
    &\int_\ell\int_\ell|t_1-t_2|^{n-1}1_E(t_1)1_E(t_2)dt_1dt_2\nonumber\\
    =&\int_0^{|E_\ell|}\int_0^{|E_\ell|}|I(t_1)-I(t_2)|^{n-1}dt_1dt_2\nonumber\\
    \ge&\int_0^{|E_\ell|}\int_0^{|E_\ell|}|t_1-t_2|^{n-1}dt_1dt_2\label{eq:ineq}\\
    =&\ \frac{1}{n}\int_0^{|E_\ell|}(t_2^n+(|E_\ell|-t_2)^n)dt_2\nonumber\\
    =&\ \frac{2}{n(n+1)}|E_\ell|^{n+1}.\label{eq:inner_int_E^2}
\end{align}
Combine \eqref{eq:E^2}\eqref{eq:inner_int_E^2}, we get the desired inequality $c(E)\le1$. Note that $c(E)\ge0$ trivially holds.

For convex $E$, we see that equality in \eqref{eq:ineq} holds for all $\ell\in\Omega$, so $c(E)=1$. More generally, if $E$ differs from a convex set on a set of measure zero, then $c(E)=1$ by Proposition~\ref{prop:index}.

Conversely, when $c(E)=1$, we see that for almost every $\ell\in\Omega$, the equality in \eqref{eq:ineq} holds. Equivalently, there exists $S\subset\S^{n-1}/\{\pm1\}$ with $|S|=|\S^{n-1}/\{\pm1\}|$, such that for all $\xi\in S$, the equality in \eqref{eq:ineq} holds for almost every $\ell$ in direction $\xi$. For these $\ell$, since $I$ is monotonically increasing with \eqref{eq:Lip}, we deduce that 
\[
    I(t_1)-I(t_2)=t_1-t_2\text{ for all }t_1,t_2\in(0,|E_\ell|).
\]
This implies that there exists a segment $A(\ell)\subset\ell$ with
\[
|A(\ell)|=|E_\ell|=|A(\ell)\cap E_\ell|.
\]

A point $x\in\R^n$ is said to be a point of Lebesgue density of $E$ if
\[
\lim_{\substack{r(B)\to0\\x\in B}}\frac{|B\cap E|}{|B|}=1,
\]
where $B$ ranges over all balls containing $x$ and $r(B)$ denotes the radius of $B$. Let $E'\subset\R^n$ be the set of all points of Lebesgue density of $E$. It is well known (e.g. see Stein \cite[Corollary~3.1.5]{stein2009real}) that $E'$ and $E$ differ by a set of measure zero. We claim that $E'$ is convex, from which the statement follows.

Suppose $x,y\in E'$, $z$ lies on the open segment $(x,y)$. Let $\epsilon>0$ be given. Find $0<\delta_0<\min(|x-z|,|y-z|)/4$ such that $|B\cap E|/|B|>1-\epsilon/2^{n+1}$ for any ball $B$ containing $x$ or $y$ with $r(B)<\delta_0$.

Now let $B$ be any ball containing $z$ with $r(B)<\delta=\delta_0/2$. Then we can find balls $B_1$, $B_2$ containing $x,y$, respectively, both with radii $2r(B)$, such that the centers of $B,B_1,B_2$ lie on a line $\ell$ with direction in $S$. By the choice of $\delta$, we see that $B,B_1,B_2$ are disjoint. Let $\Omega_0\subset\Omega$ denote the set of all lines parallel to $\ell$ that intersect $B$. We now have
\[
\int_{\Omega_0}|(B_i\backslash E)_\ell|d\ell\le|B_i\backslash E|\le\frac{\epsilon}{2^{n+1}}|B_i|=\frac{\epsilon}{2}|B|,\ i=1,2.
\]
Denote $\Omega_i=\{\,\ell\in\Omega_0\colon|(B_i\cap E)_\ell|=0\,\}$, then
\[
2r(B)|\Omega_i|\le\int_{\Omega_i}|(B_i)_\ell|d\ell=\int_{\Omega_i}|(B_i\backslash E)_\ell|d\ell\le\frac{\epsilon}{2}|B|,\ i=1,2.
\]
By the choice of $\xi$, we see that for almost every $\ell\in\Omega_0\backslash(\Omega_1\cup\Omega_2)$, $A(\ell)$ is defined and intersects with both $B_1$ and $B_2$, thus it contains $B_\ell$. Now
\[
|B\cap E|=|B|-|B\backslash E|\ge|B|-\int_{\Omega_1\cup\Omega_2}|E_\ell|\ge|B|-|\Omega_1\cup\Omega_2|\cdot2r(B)\ge(1-\epsilon)|B|.
\]
Since $\epsilon>0$ is arbitrary, we conclude that $z\in E'$, and the proof is complete.
\end{pf}

The following theorem is an asymptotic version of Theorem~\ref{thm:index}(2). Neither of these two holds for the Beer convexity index.

\begin{Thm}\label{thm:convex}
Let $E\subset\R^n$ be a set whose convexity index is defined. Suppose there is a convex set $F\subset\R^n$ with $c|F|\le|E|\le c^{-1}|E\cap F|$ for some constant $0<c<1$, then $c(E)\sim1$. Conversely, suppose $c(E)\ge c$ for some constant $0<c<1$, then there exists a convex set $F\subset\R^n$ with $|F|\sim|E|\sim|E\cap F|$. Here all implicit constants depend only on $c,n$.
\end{Thm}

\begin{pf}
The forward direction is straightforward. We first find a ellipsoid $B\subset\R^n$ containing $F$ with $|B|\sim|F|$ (see \cite{GRUBER1993319}). Since the assumptions and conclusions are all invariant under non-singular affine transformations, we may assume $B$ is a ball with radius $1$. Then $|E|\sim1$.

We have
\begin{equation}\label{eq:X1}
    \int_\Omega|E_\ell|^{n+1}d\ell\ge\int_\Omega|(E\cap F)_\ell|^{n+1}d\ell=\int_{\S^{n-1}/\{\pm1\}}\int_{B_\xi}|(E\cap F)_{\ell(\xi,x)}|^{n+1}dxd\xi,
\end{equation}
where $B_\xi$ is the projection of $B$ along $\xi$ onto $\xi^\perp$, which is an $(n-1)$-ball of radius $1$, and $\ell(\xi,x)$ is the line passing $x$ with direction $\xi$.

By H\"older's inequality,
\[
(\int_{B_\xi}|(E\cap F)_{\ell(\xi,x)}|^{n+1}dx)^\frac{1}{n+1}|B_\xi|^\frac{n}{n+1}\ge\int_{B_\xi}|(E\cap F)_{\ell(\xi,x)}|dx=|E\cap F|\sim1,
\]
thus \eqref{eq:X1} yields
\[
    \int_\Omega|E_\ell|^{n+1}d\ell\gtrsim\int_{\S^{n-1}/\{\pm1\}}1d\xi\sim1.
\]
So $c(E)\sim1$ by definition.

For the converse direction, we first prove the following lemma. We say a set $F\subset\R^n$ is star-convex at point $x$ if $F$ contains all segments $[x,y]$, $y\in F$.

\begin{Lem}\label{lem:star-convex}
Let the assumptions be as in Theorem~\ref{thm:convex}, then there exists $E'\subset E$ with $|E'|\sim|E|$, such that for every $x\in E'$, there exists a set $F_x\subset\R^n$ that is star-convex at point $x$, with $|F_x|\sim|E|\sim|F_x\cap E|$. Here the implicit constants depend only on $c,n$.
\end{Lem}

\begin{pf}
For $x\in\R^n$, $\omega\in\S^{n-1}$, $E\subset\R^n$, we use $\Omega_x$ to denote the set of all lines in $\R^n$ passing $x$ with the natural measure, and $E_{x,\omega}$ to denote the intersection of $E$ with the ray with initial point $x$ and direction $\omega$. Then
\begin{align}\label{eq:x1}
    \int_\Omega|E_\ell|^{n+1}d\ell=&\int_\Omega\int_\ell1_E(x)|E_\ell|^ndxd\ell\nonumber\\
    =&\int_{\R^n}\int_{\Omega_x}1_E(x)|E_\ell|^nd\ell dx\nonumber\\
    =&\int_E\int_{\Omega_x}|E_\ell|^nd\ell dx\nonumber\\
    \sim&\int_E\int_{\S^{n-1}}|E_{x,\omega}|^nd\omega dx.
\end{align}
On the other hand, using polar coordinate, we can write
\begin{equation}\label{eq:e2}
    |E|^2\sim\int_E\int_{\S^{n-1}}\int_0^\infty r^{n-1}1_E(x+r\omega)drd\omega dx.
\end{equation}
Now \eqref{eq:x1}\eqref{eq:e2} together with the assumption $c(E)\sim1$ yield
\[
\int_E\int_{\S^{n-1}}\int_0^\infty r^{n-1}1_E(x+r\omega)drd\omega dx\sim\int_E\int_{\S^{n-1}}|E_{x,\omega}|^nd\omega dx.
\]
But for the integrands above we have
\[
\int_{\S^{n-1}}\int_0^\infty r^{n-1}1_E(x+r\omega)drd\omega\gtrsim\int_{\S^{n-1}}|E_{x,\omega}|^nd\omega.
\]
And note that the value of left hand side above is independent of $x$. Thus for generic $x\in E$ these two terms should be comparable. More precisely, there exists $E'\subset E$, $|E'|\sim|E|$ such that 
\[
\int_{\S^{n-1}}\int_0^\infty r^{n-1}1_E(x+r\omega)drd\omega\sim\int_{\S^{n-1}}|E_{x,\omega}|^nd\omega
\]
for all $x\in E'$.
Again, for the integrands above we have
\[
\int_0^\infty r^{n-1}1_E(x+r\omega)dr\gtrsim|E_{x,\omega}|^n.
\]
This time we deduce that there exists $S_1\subset\S^{n-1}$ such that
\[
\int_{S_1}|E_{x,\omega}|^nd\omega\sim\int_{\S^{n-1}}|E_{x,\omega}|^nd\omega=|E|,
\]
and that
\begin{equation}\label{eq:13}
    \int_0^\infty r^{n-1}1_E(x+r\omega)dr\sim|E_{x,\omega}|^n
\end{equation}
for all $\omega\in S_1$.

Now, \eqref{eq:13} implies that there exists $R_{x,\omega}\sim|E_{x,\omega}|$ such that $|\{\,y\in E_{x,\omega}\colon|y-x|\le R_{x,\omega}\,\}|\sim|E_{x,\omega}|$ for all $\omega\in S_1$. Now we take
\[
F_x=\{\,x+t\omega\colon\omega\in S_1,\ 0\le t\le R_{x,\omega}\,\}.
\]
Then $F_x$ is star-convex at $x$, and we have
\[
|F_x|\sim\int_{S_1}\int_0^{R_{x,\omega}}r^{n-1}drd\omega\sim\int_{S_1}(R_{x,\omega})^nd\omega\sim\int_{S_1}|E_{x,\omega}|^nd\omega\sim|E|
\]
and that
\[
|F_x\cap E|=\int_{S_1}\int_0^{R_{x,\omega}}r^{n-1}1_E(x+r\omega)drd\omega\sim\int_{S_1}|E_{x,\omega}|^n\sim|E|.
\]
The statement follows.
\end{pf}

Now we come back to the proof of the converse direction of Theorem~\ref{thm:convex}.

Since all the assumptions and conclusions are invariant under non-singular affine transformations, we may assume that $|E|=1$.

For any set $A\subset\R^n$ and point $x\in\R^n$, we use $\star(x,A)$ to denote the union of all segments $[x,y]$, $y\in A$. Use the notations in Lemma~\ref{lem:star-convex}, for each $x\in E'$, denote $E_x=|F_x\cap E|$, then $E_x\subset\star(x,E_x)\subset F_x$, so 
\begin{equation}\label{eq:size=1}
    |E_x|\sim|\star(x,E_x)|\sim|E|=1.
\end{equation}

Apply Lemma~\ref{lem:large_simplex} with $(E,\M,\mu)$ being the restriction of the Lebesgue measure space $\R^n$ to the set $E$ and $S$ being $E'$, we see that there exists $x_0,\cdots,x_n\in E'$ with
\[
\left|\bigcap_{i=0}^nE_{x_i}\right|\sim1
\]
and
\[
\vol(x_0,\cdots,x_n)\gtrsim|E'|\sim1.
\]

Denote $E_0=\bigcap_{i=0}^nE_{x_i}$.

For any set $A\subset\R^n$ and point $x\in\R^n$, we use $\cone(x,A)$ to denote the union of all rays with initial point $x$ that passes some point $y\ne x$ in $A$. We now state the following lemma which will be useful.

\begin{Lem}\label{lem:small_cone}
Suppose $G\subset\R^n$ with $|G|\sim1$, $y_1,\cdots,y_n$ are the vertices of a regular $(n-1)$-simplex contained in $\R^{n-1}\times\{0\}$ with side length $a\gtrsim1$ and center $0$. Suppose furthermore that $|\star(y_i,G)|\sim1$ for all $1\le i\le n$. Let $c_0<1/2$, $c_0\sim1$ be a constant, $\Omega_{c_0}=\{\,\omega\in\S^{n-1}\colon|\omega\cdot e_n|\ge c_0\,\}$. Then for some $c'\sim1$ and all $R>c'a$, we have
\[
|\cone(0,\Omega_{c_0})\cap G\backslash B(0,R)|\lesssim R^{-\frac{1}{2n-1}}.
\]
\end{Lem}

\begin{Cor}\label{cor:small_ell}
Let the notations be as in Lemma~\ref{lem:small_cone}. Assume $a>2$. Let $E(R,a)$ denote the ellipsoid with center $0$ whose axes agree with the coordinate axes, and whose semi-axes have length $aR$ in the direction of the first $n-1$ coordinates and $R/a^{n-1}$ in the last. Then for some $c'\sim1$ and all $R>c'$, we have
\[
|(\cone(0,\Omega_{c_0/{a^{n-1}}})\backslash E(R,a))\cap G|\lesssim R^{-\frac{1}{2n-1}}.
\]
\end{Cor}
\begin{pf}
This follows by applying the affine transformation $(\tilde x,z)\mapsto(a\tilde x,z/a^{n-1})$ ($\tilde x\in\R^{n-1},z\in\R$) to the case $a=1$ in Lemma~\ref{lem:small_cone} and some rescaling.
\end{pf}

Assuming Lemma~\ref{lem:small_cone}, we now finish the proof of the converse direction of Theorem~\ref{thm:convex}.

By applying an affine transformation with determinant $1$, we may assume $x_0,\cdots,x_n$ are vertices of a regular $n$-simplex $\Lambda$ with side length $a\gtrsim1$ and center $0$.

First we suppose $a\ge a_0$ for a large constant $a_0>2$.

Denote the $(n-1)$-hypersurfaces of the simplex $\Lambda$ by $\alpha_0,\cdots,\alpha_n$. Then
\[
\R^n=\bigcup_{i=0}^n\cone(0,\alpha_i)=\bigcup_{i=0}^n(-\cone(0,\alpha_i)).
\]
But each $-\cone(0,\alpha_i)$ is contained in the associated ``cone substracting ellipsoid'' of $\alpha_i$ given by Corollary~\ref{cor:small_ell} if we take $G=E_0$, $c_0$ suitably small, $R=c'a$, and $a_0$ suitably large. It follows that $1\sim|E_0|\lesssim a^{-\frac{1}{2n-1}}\lesssim a_0^{-\frac{1}{2n-1}}$, so $a_0\lesssim1$.

Below we can assume $a\le a_0$ for a fixed large constant $a_0$.

Take $G=E_0$, $c_0>0$ suitably small, $R=(1+c')a$ in Lemma~\ref{lem:small_cone}, then by the same argument, the corresponding $n+1$ cones associated with $\alpha_0,\cdots,\alpha_n$ cover all of $\R^n$. It follows that
\[
|E_0\backslash B(0,2R)|\le|E_0\backslash B(0,R+a)|\lesssim R^{-\frac{1}{2n-1}}
\]
for all $R>c'a$. Hence we can find a constant $R_0>0$ such that
\[
|E_0\backslash B(0,R_0)|\le\frac{|E_0|}{2}
\]
uniformly in $E$. Take $F=B(0,R_0)$, then $|F|\sim1$, and
\[
|E\cap F|\ge|E_0\cap F|=|E_0|-|E_0\backslash B(0,R_0)|\ge\frac{|E_0|}{2}\sim1,
\]
and the proof is complete.

\vspace{15pt}

It remains to prove Lemma~\ref{lem:small_cone}.

For convenience, we replace the given $y_1,\cdots,y_n$ by $y_0=0$, $y_i=ae_i$, $1\le i\le n-1$ and prove the same result. This reduction is justified by taking an affine transformation of determinant $1$ and some rescaling.

We prove the lemma for $c'=10/c_0$.

Let $B_1=B(0,R)$, $B_2=B(0,2R)$, $S_R=\partial B_1$, $K=(B_2\backslash B_1)\cap\cone(0,\Omega_{c_0})$, $G_R=G\cap K$. We now estimate $|G_R|$.

For $i=0,\cdots,n-1$, let $S_i=\star(y_i,G_R)\cap S_R$, $S=\bigcup_{i=0}^{n-1}S_i$. Since $\star(y_i,S_i)\subset\star(y_i,G_R)\subset\star(y_i,G)$, we have
\[
R|S_i|\sim|\star(y_i,S_i)|\le|\star(y_i,G)|\sim1.
\]
So that $|S_i|\lesssim1/R$ for all $i=0,\cdots,n-1$, and thus $|S|\lesssim1/R$.

Also, for all $i=0,\cdots,n-1$, we have
\[
G_R\subset\cone(y_i,S_i)\subset\cone(y_i,S).
\]
Thus
\begin{align}\label{eq:G_R}
    G_R\subset&\left(\bigcap_{i=0}^{n-1}\cone(y_i,S)\right)\cap K\nonumber\\
    =&\left(\bigcap_{i=0}^{n-1}\left(\bigcup_{z\in S}\left[y_i,z\right>\right)\right)\cap K\nonumber\\
    =&\left(\bigcup_{z_0,\cdots,z_{n-1}\in S}\left(\bigcap_{i=0}^{n-1}\left[y_i,z_i\right>\right)\right)\cap K.
\end{align}
Here $\left[y,z\right>$ denotes the ray with initial point $y$ that passes $z$. Notice, in the last line, $\bigcap_{i=0}^{n-1}\left[y_i,z_i\right>$ is empty if $|z_i-z_j|\ge a$ for some $0\le i,j\le n-1$. Thus, $\bigcap_{i=0}^{n-1}\left[y_i,z_i\right>\ne\emptyset$ only if all $z_i$ are contained in a cap on $S_R$ of radius $a$ (use the canonical distance on $S_R$).

Now we write $S_R=\bigcup_{\alpha\in\A} S_\alpha$, where $\{S_\alpha\}_{\alpha\in\A}$ is a collection of caps of radius $2a$ on the sphere $S_R$ with $\#\A\sim R^{n-1}/a^{n-1}$ such that each cap on $S_R$ with radius $a$ is contained some $S_\alpha$. Now \eqref{eq:G_R} further gives
\begin{align}\label{eq:alpha}
    G_R\subset&\left(\bigcup_{\alpha\in\A}\left(\bigcap_{i=0}^{n-1}\cone(y_i,S_\alpha\cap S)\right)\right)\cap K\nonumber\\
    =&\bigcup_{\alpha\in\A}\left(\left(\bigcap_{i=0}^{n-1}\cone(y_i,S_\alpha\cap S)\right)\cap K\right)\nonumber\\
    \triangleq&\bigcup_{\alpha\in\A}G_\alpha.
\end{align}
Fix an $\alpha\in\A$. We now estimate $|G_\alpha|$.

Denote the unit normal vector in the direction from $0$ to the center of $S_\alpha$ by $u$. Since $G_R\subset \cone(0,\Omega_{c_0})\backslash B_1$ and $R>10a/c_0$, one can check that $S_\alpha\cap S=\emptyset$ (which implies $G_\alpha=\emptyset$) unless $|u\cdot e_n|\ge c_0/2$. Below we shall assume $|u\cdot e_n|\ge c_0/2$. By symmetry we further assume $c_1:=u\cdot e_n\ge c_0/2$.

Denote $\pi_t=\{\,x\in\R^n\colon x\cdot e_n=t\,\}$. Take $R_1=c_1R/2<c_1R-2a$, $R_2=4c_1R>2c_1R+6a$. Then one can check that $G_\alpha$ is contained in the plank $K'$ between $\pi_{R_1}$ and $\pi_{R_2}$.

Let
\[
D=\bigcup_{i=0}^{n-1}\left(\star(y_i,S_\alpha\cap S)\cap\pi_{R_1}\right).
\]\
Then $|D|\sim|S_\alpha\cap S|$, $\diam(D)\sim a$, and that
\[
G_\alpha\subset\left(\bigcap_{i=0}^{n-1}\cone(y_i,D)\right)\cap K'.
\]
Denote $G_{i,t}=\cone(y_i,D)\cap\pi_t$, $G_t=\bigcap_{i=0}^{n-1}G_{i,t}$, then
\[
G_{i,t}=y_i+\frac{t}{R_1}(D-y_i)=\frac{t}{R_1}D-\frac{t-R_1}{R_1}y_i.
\]
Thus
\[
\frac{R_1}{t}G_{0,t}=D
\]
and
\[
\frac{R_1}{t}G_{i,t}=D-a(1-\frac{R_1}{t})e_i\text{ for }1\le i\le n-1.
\]
Now we have
\begin{align*}
    |G_\alpha|\le&\int_{R_1}^{R_2}|G_t|dt\\
    \sim&\int_{R_1}^{R_2}\left|\frac{R_1}{t}G_t\right|dt\\
    =&\int_{R_1}^{R_2}\left|D\cap\left(\bigcap_{i=1}^{n-1}\left(D-a\left(1-\frac{R_1}{t}\right)e_i\right)\right)\right|dt\\
    =&\int_0^{\frac{a(R_2-R_1)}{R_2}}\left|D\cap\left(\bigcap_{i=1}^{n-1}(D-ue_i)\right)\right|\frac{aR_1}{(a-u)^2}du\\
    \sim&\ \frac{R}{a}\int_0^{\frac{a(R_2-R_1)}{R_2}}\left|D\cap\left(\bigcap_{i=1}^{n-1}(D-ue_i)\right)\right|du\\
    \lesssim&\ \frac{R}{a}(a^{n-1}|D|^{2n})^\frac{1}{2n-1}\text{ (by Lemma~\ref{lem:small_int})}\\
    =&\ Ra^{-\frac{n}{2n-1}}|D|^\frac{2n}{2n-1}\\
    \lesssim&\ R|S_\alpha\cap S|^{\frac{2n}{2n-1}}\\
    \lesssim&\ R^{\frac{2n-2}{2n-1}}|S_\alpha\cap S|.
\end{align*}
Summing up all $\alpha$ in \eqref{eq:alpha}, we get
\[
|G_R|\lesssim R^{\frac{2n-2}{2n-1}}\sum_{\alpha\in\A}|S_\alpha\cap S|\sim R^\frac{2n-2}{2n-1}|S|\lesssim R^{-\frac{1}{2n-1}}.
\]
Thus
\[
|\cone(0,\Omega_{c_0})\cap G\backslash B(0,R)|=\sum_{k=0}^\infty|G_{2^kR}|\lesssim\sum_{k=0}^\infty(2^kR)^{-\frac{1}{2n-1}}\lesssim R^{-\frac{1}{2n-1}}.
\]
The proof for Lemma~\ref{lem:small_cone}, and thus for Theorem~\ref{thm:convex}, is complete.
\end{pf}

\section{Proof of Theorem 4.5}\label{pf}
We can first find some balls $B_1,\cdots,B_k$ of radius $1$ such that $3B_i$ are pairwise disjoint and that
\[
\sum_{i=1}^k\int_{B_i}\sum_{T\in\T}1_T\sim N\delta^{n-1}.
\]
Denote $\T_i=\{\,T\in\T\colon T\subset3B_i\,\}$, $N_i=\#\T_i$. Then we have $\sum_{i=1}^kN_i\sim N$, and that by \eqref{eq:main}\eqref{eq:small} (here and below, the implicit constant may depend on $C$),
\[
\sum_{i=1}^k\sqrt{N_i}\delta^{n-1}\lesssim\sum_{i=1}^k\left|\bigcup_{T\in\T_i}T\right|\le\left|\bigcup_{T\in\T}T\right|\lesssim\sqrt{N}\delta^{n-1}.
\]
Assume $N_1=\max_{1\le i\le k}N_i$. Then
\[
\sqrt{N}\gtrsim\sum_{i=1}^k\sqrt{N_i}\ge\sum_{i=1}^k\frac{N_i}{\sqrt{N_1}}\sim\frac{N}{\sqrt{N_1}}.
\]
Thus $N_1\sim N$, and $|\bigcup_{T\in\T_1}T|\le|\bigcup_{T\in\T}T|\lesssim\sqrt{N}\delta^{n-1}\sim\sqrt{N_1}\delta^{n-1}$. Now, we see that the theorem follows if we can prove the statement for $\T$ replaced by $\T_1$ (with a possibly larger constant $C$).

Thus, we may as well assume at the beginning that all $T\in\T$ are contained in a ball of radius $3$.

Then, by an application of pigeonhole principle, we may assume that all tubes have direction contained in a $1/100$-cap. Moreover, using the trick in the proof of Lemma~\ref{lem:center}, we may also assume there is a point $O$ such that all tubes are within $1/100$ vertical distance with respect to $O$. Without loss of generality let $O$ be the origin, assume all tubes have direction within $1/100$ of the vertical direction (the last coordinate).

At this point, we use the arithmetic method first introduced by Bourgain \cite{bourgain1999dimension} to prove the following result.
\begin{Prop}
Suppose \eqref{eq:small} holds. Then there exists a collection $\I$ of $\delta$-separated directions with $\#\I\sim\sqrt{N}$ such that $\#\T_e\sim\sqrt{N}$ for all $e\in\I$, where 
\[
\#\T_e=\{\,T\in\T\colon T \text{ has direction in }\c(e,\frac{\delta}{2})\,\}.
\]
\end{Prop}
\begin{pf}
We may assume $\T$ has the properties as remarked above.

Denote $K=\bigcup_{T\in\T}T$, $\pi_t=\{\,x\in\R^n\colon x\cdot e_n=t\,\}$. Then for all $|t|\le1/4$, $\pi_t$ intersect each tube in $\T$ at an $(n-1)$-ellipsoid whose axes all have length within $[\delta,2\delta]$. Let $K_t=K_{\pi_t}$, $X=\{\,t\in\R\colon|K_t|>100|=|K|\,\}$, then $|X|\le1/100$.

Pick any $t_0\in[-1/4,0]\backslash X$, then
\[
\int_{1/16}^{1/8}(1_X(t_0+d_0)+1_X(t_0+2d_0))dd_0\le|X|\le1/100
\]
tells us that there exists $d_0\in[1/16,1/8]$ such that $t_0+d_0,t_0+2d_0\not\in X$. Denote $t_i=t_0+id_0$. Then $|t_i|\le1/4$ and $|K_{t_i}|\lesssim|K|$, $i=0,1,2$.

For $|t|\le1/4$, let $Z_t'=c_0\Z^{n-1}\times\{t\}\subset\pi_t$, $Z_t=2c_0\Z^{n-1}\times\{t\}\subset Z_t'$. Then for a suitably chosen small constant $c_0\sim1$, we have 
\[
\#(K_t\cap Z_t)\sim\#(K_t\cap Z_t')\sim|K_t|\delta^{1-n}.
\]

Now pick any $a(T)\in T\cap Z_{t_0}$, $b(T)\in T\cap Z_{t_2}$ for each $T\in\T$. Then $(a(T)+b(T))/2\in T\cap Z_{t_1}'$. Moreover, for any $(a,b)\in Z_{t_0}\times Z_{t_2}$, by pigeonhole principle (note that $d_0\sim1$) there are $\lesssim1$ tubes $T\in\T$ with $a(T)=a$, $b(T)=b$.

Identify $Z_t$ with its projection in $Z=2c_0\Z^{n-1}$. Let $A=\{\,a(T)\colon T\in\T\,\}$, $B=\{\,b(T)\colon T\in\T\,\}$, $C=\{\,a(T)+b(T)\colon T\in\T\,\}$, $G=\{\,(a(T),b(T))\colon T\in\T\,\}$. Then by the remarks above we have
\[
\#A\le\#(K_{t_0}\cap Z_{t_0})\sim|K_{t_0}|\delta^{1-n}\lesssim|K|\delta^{1-n}\lesssim N^{1/2},\ \#B\lesssim N^{1/2},
\]
\[
\#C\le\#(K_{t_1}\cap Z_{t_1}')\lesssim N^{1/2},\ \#G\sim\#\T=N.
\]
Suppose $T\mapsto\varphi(T)=a(T)-b(T)$ maps at most $m$ points to one point. Apply Lemma~\ref{lem:arithmetic}, we see that
\[
N\lesssim m^\frac{1}{6}(N^\frac{1}{2})^\frac{11}{6}.
\]
In other words $m\gtrsim\sqrt{N}$. Thus there exists some $u\in Z$ being the image of $\gtrsim\sqrt{N}$ tubes $T\in\T$ under $\varphi$. This means there are $\gtrsim\sqrt{N}$ tubes with direction within $\lesssim\delta$ of the direction determined by the vector $(u,2d_0)$ (note again that $d_0\sim1$). By pigeonhole principle there exists a direction $e\in\S^{n-1}/\{\pm1\}$ with $\#\T_e\gtrsim\sqrt{N}$.

On the other hand, by definition of $\T_e$ and pigeonhole principle,
\[
\#\T_e\delta^{n-1}\sim\left|\bigcup_{T\in\T_e}T\right|\le|K|\lesssim\sqrt{N}\delta^{n-1},
\]
so $\#\T_e\sim\sqrt{N}$.

Above, all implicit constants depend only on $n,C$. Now if we enlarge $C$ to $2C$ and do the same argument, we see that for any subcollection $\T'$ of $\T$ with cardinality $\ge N/4$, we can always find a direction $e\in\S^{n-1}/\{\pm1\}$ with $\T'_e\sim\sqrt{N}$. Hence the proposition follows when we repeat this argument $\sim\sqrt{N}$ times.
\end{pf}

Now back to the proof of Theorem~\ref{thm:small}. The theorem follows if we prove the corresponding result for $\T$ replaced by $\bigcup_{e\in\I}\T_e$, which we now assume. Denote $K=\bigcup_{T\in\T}T$. Let $\mu=\sum_{T\in\T}1_T$, then $\mu\lesssim\sqrt{N}$ everywhere by the bush argument in the proof of Theorem~\ref{thm:bound}. We have
\[
N\delta^{n-1}\sim|K|\sqrt{N}\gtrsim\int_K\mu=\sum_{T\in\T}\int_{\R^n}1_T\sim N\delta^{n-1}.
\]
Equality implies $\mu\sim\sqrt{N}$ in some $K'\subset K$ with $|K'|\sim|K|$. Thus
\[
N^\frac{3}{2}\delta^{n-1}\sim\int_K\mu^2=\sum_{T\in\T}\int_T\mu\lesssim\#\T\sqrt{N}\delta^{n-1}=N^\frac{3}{2}\delta^{n-1}.
\]
Again, equality implies that there is some $\T_0\subset\T$ with $\#\T_0\sim N$ such that
\begin{equation}\label{eq:T_0}
    \int_T\mu\sim\sqrt{N}\delta^{n-1}\text{ for all }T\in\T_0.
\end{equation}
Let $\T_e'=\T_e\cap\T_0$, $e\in\I$. We can choose a certain $e\in\I$ such that $\#\T_e'\sim\sqrt{N}$.

By applying a rotation (and loosen the constant $1/100$ to $1/50$ in the assumption that all tubes have direction within $1/100$ of the vertical direction) if necessary, we assume $e$ is the vertical direction. For $T\in\T_e'$ with some center $a$, let $\widetilde T$ be the tube with vertical direction, radius $2\delta$, height $2$, and center $\Pi_n(a)$, where $\Pi_n$ is the projection map onto $\R^{n-1}\times\{0\}$. Recall $T$ is within vertical distance $1/100$ with respect to $0$, we see that $T\subset\widetilde T$. Denote $E_0=\Pi_n(\bigcup_{T\in\T'_e}\widetilde T)$ identified with a subset of $\R^{n-1}$. Then $\bigcup_{T\in\T'_e}\widetilde T=E_0\times[-1,1]$, $E_0$ is $2\delta$-discretized, $\diam(E_0)\le10$, $|E_0|\sim\sqrt{N}\delta^{n-1}$, and that by \eqref{eq:T_0},
\begin{equation}\label{eq:TT}
    \sum_{T\in\T}\int_{E_0\times[-1,1]}1_T=\int_{E_0\times[-1,1]}\mu\ge\int_{\bigcup_{T\in\T'_e}T}\mu\sim N\delta^{n-1}.
\end{equation}
On the other hand,
\[
\int_{E_0\times[-1,1]}1_T=|T\cap(E_0\times[-1,1])|\le|T|\sim\delta^{n-1}
\]
holds for all $T\in\T$. Now equality in \eqref{eq:TT} implies there exists some $\T'\subset\T$ with $\#\T'\sim N$ such that 
\begin{equation}\label{eq:T'}
    |T\cap(E_0\times[-1,1])|\sim|T|\text{ for all }T\in\T'.    
\end{equation}

For any set $A\subset\R$, $0<\lambda<1$, denote 
\[m(A,\lambda)=\sup\{\,|I|\colon I\subset\R\text{ is an interval, }|I\cap A|\ge\lambda|I|\,\}.
\]
Under this notation, a similar continuization argument (in the reverse direction) as in the proof of Theorem~\ref{thm:convex_is_good} shows that
\begin{equation}\label{eq:int_m}
    N\sim\#\T'\lesssim\delta^{2-2n}\int_{\Omega_{n-1}}m((E_0)_\ell,\lambda)^nd\ell,
\end{equation}
where $\lambda$ is a constant depends only on $n$ and the implicit constant in \eqref{eq:T'}, that is to say, $\lambda\sim1$.

Next, we show that there exists convex $F\subset\R^{n-1}$ with $|F|\sim|F\cap E_0|\sim\sqrt{N}\delta^{n-1}$.

\textit{Case 1:} $n=2$.

Then \eqref{eq:int_m} becomes $m(E_0,\lambda)\gtrsim\sqrt{N}\delta$. Thus we can find an interval $F\subset\R$ such that $|F|\sim|F\cap E_0|\gtrsim\sqrt{N}\delta$ by definition of $m$. This is actually an (approximate) equality since $|F\cap E_0|\le|E_0|\sim\sqrt{N}\delta$.

\textit{Case 2:} $n\ge3$.

By definition of $m$, we see that $m(A,\lambda)\le|A|/\lambda$. Thus \eqref{eq:int_m} further gives
\[
N\lesssim\delta^{2-2n}\int_{\Omega_{n-1}}|(E_0)_\ell|^n.
\]
This means $c_{n-1}(E_0)\sim1$. Now the existence of $F$ follows from Theorem~\ref{thm:convex}.

\vspace{15pt}
Let us finally finish the proof of Theorem~\ref{thm:small}.

By replacing $F$ by a slightly larger set if necessary, we may assume that $F$ is a box in $\R^{n-1}$ (see \cite{GRUBER1993319}). Suppose the shortest side of $F$ has length $d$. Then the intersection of any ball of radius $2\delta$ in $\R^{n-1}$ with $F$ has volume $\lesssim d\delta^{n-2}$. Since $E_0$ is the union of $\sim\sqrt{N}$ such balls, we see that $|E_0\cap F|\lesssim\sqrt{N}d\delta^{n-2}$. By choice of $F$ it follows that $d\gtrsim\delta$. 

Since $|(F\times[-1,1])\cap(\bigcup_{T\in\T'_e}\widetilde T)|\sim|F\cap E_0|\sim\sqrt{N}\delta^{n-1}$, we see that $F\times[-1,1]$ intersects with $\sim\sqrt{N}$ tubes $\widetilde T$ in the union. Replacing $F$ by a slightly larger box if necessary (notice $d\gtrsim\delta$, so this will not change $|F|$ too much), we may assume $F\times[-1,1]$ contains $\sim\sqrt{N}$ such tubes $\widetilde T$, and thus contains $\sim\sqrt{N}$ tubes in $\T'_e$.

Since $\T'_e\subset\T_0$, by \eqref{eq:T_0} we see that
\begin{equation}\label{eq:F_1}
    \sum_{T\in\T}\int_{F\times[-1,1]}1_T=\int_{F\times[-1,1]}\mu\gtrsim\sum_{\substack{T\in\T'_e\\T\subset F\times[-1,1]}}\int_T\mu\gtrsim N\delta^{n-1}.
\end{equation}
But $\int_{F\times[-1,1]}1_T\lesssim\delta^{n-1}$ for all $T\in\T$. Equality in \eqref{eq:F_1} shows that there exists $\lambda_0\sim1$ and $\T''\subset\T$ with $\#\T''\sim N$, such that
\[
|T\cap(F\times[-1,1])|\ge\lambda_0|T|\text{ for all }T\in\T''.
\]
For all $T\in\T''$, this further implies there is a segment, say $I$, of length $1$ contained in $T$ with the same direction as $T$, such that $|I\cap(F\times[-1,1])|\ge\lambda_0$. Then $\Pi_n(I)$ is easily seen to be contained in the box
\[
F_1=\frac{1}{\lambda_0}F+\frac{\lambda_0-1}{\lambda_0}F,
\]
where ``$+$'' denotes the Minkowski sum of sets. Hence, $T$ is contained in $F_2\times[-1,1]$, where $T_2$ is the $\delta$-neighborhood of $F_1$.

Now $|F_2|\sim|F_1|\sim|F|\sim\sqrt{N}\delta^{n-1}$ and $F_2\times[-1,1]$ contains all tubes in $\T''$, where $\#\T''\sim N$. Since all $T\in\T$ are contained in a ball of radius $3$ and have direction within $1/50$ of the vertical direction, we can find a box $E'\subset F_2$ with $|E'|\sim\sqrt{N}\delta^{n-1}$ and $\diam(E')<1/2$ such that $E'\times[-1,1]$ contains $\sim N$ tubes in $\T$. Let $E$ be the $9\delta$-neighborhood of $E'$, then $E$ satisfies all the requirements in our statement. The proof is complete.\qed

\section{Acknowledgements}
This research was done during the 2019 Summer Program of Undergraduate Research (SPUR) of the MIT Mathematics Department. First of all, I want to thank my mentor Yuqiu Fu who has provided many insightful ideas and spotted numerous errors I have made. I also want to thank Professor Larry Guth for the proposal of the original problem and Professor David Jerison for the suggestion of considering the inverse problem. Finally, I thank my mother for her constant support, without which I would undoubtedly fail to complete this work.

\bibliographystyle{amsplain}
\bibliography{citations}

\providecommand{\bysame}{\leavevmode\hbox to3em{\hrulefill}\thinspace}
\providecommand{\MR}{\relax\ifhmode\unskip\space\fi MR }
\providecommand{\MRhref}[2]{%
  \href{http://www.ams.org/mathscinet-getitem?mr=#1}{#2}
}
\providecommand{\href}[2]{#2}
\begin{thebibliography}{10}

\bibitem{beer1974index}
Gerald Beer, \emph{The index of convexity and parallel bodies}, Pacific Journal
  of Mathematics \textbf{53} (1974), no.~2, 337--345.

\bibitem{bourga1991besicovitch}
Jean Bourgain, \emph{Besicovitch type maximal operators and applications to
  fourier analysis}, Geometric and Functional analysis \textbf{1} (1991),
  no.~2, 147--187.

\bibitem{bourgain1999dimension}
\bysame, \emph{On the dimension of kakeya sets and related maximal
  inequalities}, Geometric and Functional Analysis \textbf{9} (1999), no.~2,
  256--282.

\bibitem{ChristXray}
Michael Christ, \emph{Estimates for the k-plane transform}, Indiana University
  Mathematics Journal \textbf{33} (1984), no.~6, 891--910.

\bibitem{GRUBER1993319}
Peter~M. Gruber, \emph{Chapter 1.10 - aspects of approximation of convex
  bodies}, Handbook of Convex Geometry (P.M. Gruber and J.M. Wills, eds.),
  North-Holland, Amsterdam, 1993, pp.~319 -- 345.

\bibitem{katz1999bounds}
Nets~Hawk Katz and Terence Tao, \emph{Bounds on arithmetic projections, and
  applications to the kakeya conjecture}, Mathematical Research Letters
  \textbf{6} (1999), no.~6, 625--630.

\bibitem{katz2002new}
\bysame, \emph{New bounds for kakeya problems}, Journal d'Analyse
  Math{\'e}matique \textbf{87} (2002), no.~1, 231--263.

\bibitem{MANILEVITSKA199319}
Peter Mani-Levitska, \emph{Chapter 1.1 - characterizations of convex sets},
  Handbook of Convex Geometry (P.M. Gruber and J.M. Wills, eds.),
  North-Holland, Amsterdam, 1993, pp.~19 -- 41.

\bibitem{Oberlin2010}
Richard Oberlin, \emph{Two bounds for the x-ray transform}, Mathematische
  Zeitschrift \textbf{266} (2010), no.~3, 623--644.

\bibitem{ren1994topics}
D.~Ren, \emph{Topics in integral geometry}, Pure Mathematics, World Scientific,
  1994.

\bibitem{stein2009real}
Elias~M Stein and Rami Shakarchi, \emph{Real analysis: measure theory,
  integration, and hilbert spaces}, Princeton University Press, 2009.

\bibitem{wolff1995improved}
Thomas~H Wolff, \emph{An improved bound for kakeya type maximal functions},
  Revista Matem{\'a}tica Iberoamericana \textbf{11} (1995), no.~3, 651--674.

\end{thebibliography}

\parskip=0pt

\end{document}